\documentclass[11pt]{article}

\usepackage[a4paper, margin=0.8in]{geometry}

\usepackage{amsmath,amsthm}
\usepackage{caption}
\usepackage{empheq}
\usepackage[ruled,lined,linesnumbered]{algorithm2e}
\usepackage[colorlinks=true,linkcolor=blue!60!black,citecolor=blue!60!black,urlcolor=blue!60!black]{hyperref}
\usepackage[capitalise,noabbrev]{cleveref}

\usepackage{amssymb}
\usepackage{mathtools}
\usepackage{mathrsfs}
\usepackage{bbm}
\let\origmathbb\mathbb
\renewcommand{\mathbb}[1]{\ifnum\pdfstrcmp{#1}{1}=0 \mathbbm{1}\else\origmathbb{#1}\fi}
\newcommand{\Lip}{\mathsf{Lip}}
\numberwithin{equation}{section}

\usepackage{enumitem}
\usepackage{booktabs}
\usepackage{graphicx}
\usepackage{float}
\usepackage{tikz}
\usetikzlibrary{positioning, arrows.meta, calc}

\usepackage[backend=biber,style=numeric,sorting=nty,doi=false,isbn=false,url=false,date=year]{biblatex}
\DeclareNameAlias{default}{family-given}
\DeclareFieldFormat[article]{title}{#1}
\renewbibmacro{in:}{}
\AtEveryBibitem{\ifentrytype{article}{\clearfield{number}\clearfield{issue}\clearfield{eprint}}{}}
\AtEveryBibitem{\ifentrytype{book}{\clearfield{series}\clearfield{volume}\clearfield{number}}{}}
\AtEveryBibitem{\ifentrytype{software}{\clearfield{note}}{}}
\ExplSyntaxOn
\seq_new:N \l__mfgbib_given_seq
\seq_new:N \l__mfgbib_giveni_seq
\tl_new:N \l__mfgbib_given_tl
\tl_new:N \l__mfgbib_out_tl
\bool_new:N \l__mfgbib_spelled_bool
\bool_new:N \l__mfgbib_previnit_bool
\prg_new_conditional:Npnn \__mfgbib_if_initial:n #1 { TF }
  {
    \bool_lazy_or:nnTF
      { \int_compare_p:nNn { \tl_count:n {#1} } = { 1 } }
      {
        \bool_lazy_and_p:nn
          { \int_compare_p:nNn { \tl_count:n {#1} } = { 2 } }
          { \str_if_eq_p:ee { \tl_item:nn {#1} { 2 } } { . } }
      }
      { \prg_return_true: } { \prg_return_false: }
  }
\cs_new_protected:Npn \__mfgbib_build_given:
  {
    \tl_clear:N \l__mfgbib_out_tl
    \tl_set_eq:NN \l__mfgbib_given_tl \namepartgiven
    \tl_if_blank:VF \l__mfgbib_given_tl
      {
        \tl_replace_all:Nnn \l__mfgbib_given_tl { \bibnamedelima } { \q_mark }
        \tl_replace_all:Nnn \l__mfgbib_given_tl { \bibnamedelimb } { \q_mark }
        \tl_replace_all:Nnn \l__mfgbib_given_tl { \bibnamedelimi } { \q_mark }
        \seq_set_split:NnV \l__mfgbib_given_seq { \q_mark } \l__mfgbib_given_tl
        \seq_set_split:NnV \l__mfgbib_giveni_seq { \bibinitdelim } \namepartgiveni
        \bool_set_false:N \l__mfgbib_spelled_bool
        \bool_set_false:N \l__mfgbib_previnit_bool
        \seq_map_indexed_inline:Nn \l__mfgbib_given_seq
          {
            \int_compare:nNnT {##1} > { 1 }
              {
                \bool_if:NTF \l__mfgbib_previnit_bool
                  { \tl_put_right:Nn \l__mfgbib_out_tl { \bibnamedelimi } }
                  {
                    \int_compare:nNnTF {##1} = { \seq_count:N \l__mfgbib_given_seq }
                      { \tl_put_right:Nn \l__mfgbib_out_tl { \bibnamedelima } }
                      { \tl_put_right:Nn \l__mfgbib_out_tl { \bibnamedelimb } }
                  }
              }
            \bool_if:NTF \l__mfgbib_spelled_bool
              {
                \tl_put_right:Ne \l__mfgbib_out_tl
                  { \seq_item:Nn \l__mfgbib_giveni_seq {##1} }
                \bool_set_true:N \l__mfgbib_previnit_bool
              }
              {
                \tl_put_right:Nn \l__mfgbib_out_tl {##2}
                \__mfgbib_if_initial:nTF {##2}
                  { \bool_set_true:N \l__mfgbib_previnit_bool }
                  {
                    \bool_set_true:N \l__mfgbib_spelled_bool
                    \bool_set_false:N \l__mfgbib_previnit_bool
                  }
              }
          }
      }
  }
\DeclareNameFormat{family-given}
  {
    \ifgiveninits
      {
        \usebibmacro{name:family-given}
          {\namepartfamily} {\namepartgiveni} {\namepartprefix} {\namepartsuffix}
      }
      {
        \__mfgbib_build_given:
        \usebibmacro{name:family-given}
          {\namepartfamily} {\l__mfgbib_out_tl} {\namepartprefix} {\namepartsuffix}
      }
    \usebibmacro{name:andothers}
  }
\ExplSyntaxOff

\usepackage{xcolor}
\definecolor{oceanicAbyss1}{HTML}{112C46}
\definecolor{oceanicCurrent1}{HTML}{274762}
\definecolor{oceanicMist1}{HTML}{8FAABE}
\definecolor{oceanicCyan1}{HTML}{4A8B99}
\definecolor{oceanicTeal1}{HTML}{4A7A7A}
\definecolor{oceanicSage1}{HTML}{5A8A6A}
\definecolor{oceanicMoss1}{HTML}{6A9A5A}
\definecolor{oceanicOlive1}{HTML}{7A8A5A}
\definecolor{oceanicCitrus1}{HTML}{ADA85C}
\definecolor{oceanicAmber1}{HTML}{E5C27D}
\definecolor{oceanicCoral1}{HTML}{C4866B}
\definecolor{oceanicRose1}{HTML}{A56B7A}
\definecolor{oceanicCrimson1}{HTML}{9A4A4A}
\definecolor{oceanicMagenta1}{HTML}{8A5A7A}
\definecolor{oceanicPlum1}{HTML}{6B5A7A}
\definecolor{oceanicViolet1}{HTML}{7B6B9D}
\definecolor{oceanicFoam1}{HTML}{C8D6D6}
\definecolor{oceanicCloud1}{HTML}{E8EEF0}
\definecolor{oceanicSand1}{HTML}{E5DCC8}
\definecolor{oceanicSlate1}{HTML}{3A4555}
\definecolor{oceanicInk1}{HTML}{1A2530}
\definecolor{oceanicAbyss2}{HTML}{002A5B}
\definecolor{oceanicCurrent2}{HTML}{004778}
\definecolor{oceanicMist2}{HTML}{80ACCC}
\definecolor{oceanicCyan2}{HTML}{0090A8}
\definecolor{oceanicTeal2}{HTML}{217F80}
\definecolor{oceanicSage2}{HTML}{3B915E}
\definecolor{oceanicMoss2}{HTML}{599F3F}
\definecolor{oceanicOlive2}{HTML}{758E3A}
\definecolor{oceanicCitrus2}{HTML}{B3A80F}
\definecolor{oceanicAmber2}{HTML}{F6BD39}
\definecolor{oceanicCoral2}{HTML}{DB784B}
\definecolor{oceanicRose2}{HTML}{B95C77}
\definecolor{oceanicCrimson2}{HTML}{AB383E}
\definecolor{oceanicMagenta2}{HTML}{994C82}
\definecolor{oceanicPlum2}{HTML}{705488}
\definecolor{oceanicViolet2}{HTML}{7F63B4}
\definecolor{oceanicFoam2}{HTML}{BCDCD4}
\definecolor{oceanicCloud2}{HTML}{DFF2F8}
\definecolor{oceanicSand2}{HTML}{F0DCAA}
\definecolor{oceanicSlate2}{HTML}{324562}
\definecolor{oceanicInk2}{HTML}{102539}
\definecolor{oceanicAbyss3}{HTML}{4171A2}
\definecolor{oceanicCurrent3}{HTML}{5791C1}
\definecolor{oceanicMist3}{HTML}{C7EAFF}
\definecolor{oceanicCyan3}{HTML}{6BCEE3}
\definecolor{oceanicTeal3}{HTML}{72BBBB}
\definecolor{oceanicSage3}{HTML}{85CE9E}
\definecolor{oceanicMoss3}{HTML}{A2DC8F}
\definecolor{oceanicOlive3}{HTML}{B3CB83}
\definecolor{oceanicCitrus3}{HTML}{F0E987}
\definecolor{oceanicAmber3}{HTML}{FFE37F}
\definecolor{oceanicCoral3}{HTML}{FFC3A6}
\definecolor{oceanicRose3}{HTML}{F69EB5}
\definecolor{oceanicCrimson3}{HTML}{E28180}
\definecolor{oceanicMagenta3}{HTML}{D58DBD}
\definecolor{oceanicPlum3}{HTML}{AB91C2}
\definecolor{oceanicViolet3}{HTML}{BBA2EE}
\definecolor{oceanicFoam3}{HTML}{D6E7E7}
\definecolor{oceanicCloud3}{HTML}{ECF3F6}
\definecolor{oceanicSand3}{HTML}{F7ECD4}
\definecolor{oceanicSlate3}{HTML}{7A8DA8}
\definecolor{oceanicInk3}{HTML}{54687B}

\DeclarePairedDelimiter{\parens}{(}{)}
\DeclarePairedDelimiter{\bracks}{[}{]}
\DeclarePairedDelimiter{\abs}{\lvert}{\rvert}
\DeclarePairedDelimiter{\norm}{\lVert}{\rVert}

\DeclarePairedDelimiterX\ensemble[1]\{\}{#1}
\newcommand{\diver}{\operatorname{div}}
\newcommand{\Id}{\mathrm{Id}}
\DeclareMathOperator*{\argmin}{arg\,min}\DeclareMathOperator*{\argmax}{arg\,max}

\newcommand{\diff}{\mathop{}\!\mathrm{d}}

\providecommand{\mbbE}{\mathbb{E}}

\providecommand{\mbbP}{\mathbb{P}}\providecommand{\mbbR}{\mathbb{R}}

\providecommand{\mcalA}{\mathcal{A}}
\providecommand{\mcalI}{\mathcal{I}}\providecommand{\mcalK}{\mathcal{K}}\providecommand{\mcalL}{\mathcal{L}}
\providecommand{\mcalN}{\mathcal{N}}\providecommand{\mcalR}{\mathcal{R}}
\providecommand{\mcalV}{\mathcal{V}}\providecommand{\mcalX}{\mathcal{X}}

\providecommand{\mscrG}{\mathscr{G}}\providecommand{\mscrI}{\mathscr{I}}\providecommand{\mscrJ}{\mathscr{J}}

\providecommand{\mscrS}{\mathscr{S}}

\providecommand{\msfP}{\mathsf{P}}

\providecommand{\mbfv}{\mathbf{v}}

\theoremstyle{plain}
\newtheorem{theorem}{Theorem}[section]
\newtheorem{proposition}[theorem]{Proposition}
\newtheorem{lemma}[theorem]{Lemma}
\newtheorem{corollary}[theorem]{Corollary}
\theoremstyle{definition}
\newtheorem{definition}[theorem]{Definition}
\theoremstyle{remark}
\newtheorem{remark}[theorem]{Remark}

\title{A Monotone Meshfree--Particle Method for Mean Field Games on Scattered Point Clouds}
\author{Jiongyi Wang$^{1}$\footnote{Email: jiongyiwang@keio.jp} \, and Linyu Peng$^{2}$\footnote{Email: l.peng@mech.keio.ac.jp}  \vspace{0.4cm}\\
{\it 1. Graduate School of Science and Technology, Keio University,} \\
{\it Yokohama 223-8522, Japan}\\
{\it 2. Department of Mechanical Engineering, Keio University,} \\
{\it Yokohama 223-8522, Japan}\\
}

\SetAlCapFnt{\footnotesize}\SetAlCapNameFnt{\footnotesize}
\begin{document}

\maketitle

\begin{abstract}
    We develop a monotone meshfree--particle method for mean field games (MFGs). We discretize the backward Hamilton--Jacobi--Bellman (HJB) equation on a fixed scattered cloud, transport the forward Fokker--Planck law with reflected and killed particles, and close the fictitious-play loop by interpolation and mass-preserving kernel density estimation. Monotonicity is imposed row by row on two-dimensional clouds. Each row combines a positive generalized finite difference Laplacian with a nonnegative decomposition of the drift along rays to surrounding visible neighbors, with mirror images at reflecting walls, and every fixed-policy matrix is an M-matrix for arbitrary cell-P\'eclet number. We prove a discrete comparison principle and, on polygons with re-entrant corners included, that one HJB solve at a prescribed density converges to the viscosity solution off the boundary vertices, for convex Hamiltonians of at most quadratic growth under a structure condition. On a two-room floor plan with re-entrant door jambs, the computed HJB solutions converge to exact ones. A source-free exact MFG tests the fully coupled computation. At matched degrees of freedom on a grid domain, the meshfree HJB solve has a larger $L^2$ error and a smaller $L^\infty$ error than a five-point finite-difference solve, and is about twice as expensive including the stencil construction. On a non-convex evacuation game, the interior nodes of nine clouds receive monotone rows, apart from a few next to the boundary vertices that copy a neighbor.

\end{abstract}

\tableofcontents


\section{Introduction}

The self-organized motion of large interacting populations, such as pedestrian crowds evacuating a floor plan or agents steering around obstacles toward a target, can be modeled, in the continuum limit, by a coupled system of a backward Hamilton--Jacobi--Bellman (HJB) equation and a forward Fokker--Planck (FP) transport equation. Mean field game (MFG) theory, introduced independently by Lasry and Lions \cite{Lasry2007Mean} and Huang, Malham\'e, and Caines \cite{Caines2006Large}, furnishes this structure as the mean-field limit of $N$-player differential games. For the pedestrian setting, it continues a modeling line running from social-force and continuum crowd models \cite{HelbingFarkasVicsek2000,Hughes2002Continuum} to their MFG formulations \cite{AchdouLasry2019Crowd}.

Among MFG numerical methods, established HJB stability and, when adjoint-consistent, Eulerian mass conservation come with a structured or body-fitted mesh. For coupled second-order MFGs with nondifferentiable Hamiltonians, convergence rates follow from a single stabilized bilinear form serving both equations \cite{Osborne2026Rates}. There, the forward operator is the discrete adjoint of the backward one, so both inherit one stability constant. The framework is posed on polytopal domains with homogeneous Dirichlet data and acute meshes. Particle approximation has been carried out for one-dimensional first-order planning problems \cite{DiFrancescoEtAl2022}, while trajectory-based semi-Lagrangian schemes reach degenerate second-order ones \cite{Carlini2015Semi2nd}. Learning and flow-matching routes bypass the classical HJB solver altogether \cite{Lauriere2022Learning,Yu2025FlowMatching}.

A monotone, stable and consistent scheme converges to the viscosity solution whenever the limit problem admits a comparison principle \cite{BarlesSouganidis1991}. On a grid or a mesh, monotonicity is a property of the whole discretization. The upwind five-point stencil \cite{Achdou2010Mean} needs a structured grid, and edge-averaged elements give M-matrices for every continuous diffusion and drift exactly when the triangulation satisfies a discrete angle condition, which, in the plane, is equivalent to the Delaunay condition up to rare exceptions near the boundary \cite[Lems.~2.1 and~3.2]{XuZikatanov1999}. One violating edge therefore voids the guarantee for the entire matrix. On a grid, wide stencils give monotone HJB diffusion terms, with nonnegative coefficients from a linear program at each node \cite{BonnansOttenwaelterZidani2004}. Body-fitted finite-element and discontinuous-Galerkin methods conform to complex domains \cite{SmearsSuli2014,Osborne2025Finite}. On the boundary side, monotone second-order schemes reach oblique conditions on smooth bounded domains \cite{CalzolaCarliniDupuisSilva2023} and nonlinear mixed conditions on unstructured meshes \cite{JaroszkowskiJensen2024}, the latter assuming comparison for the limit problem outright.

For monotone meshless discretization, Seibold \cite{Seibold2008Minimal} constructed positive Laplacian stencils, Ye and Tian \cite{YeTian2023} treated linear non-divergence-form elliptic operators, Froese \cite{Froese2018Meshfree} obtained monotone point-cloud approximations of functions of the Hessian eigenvalues converging to viscosity solutions, and Finlay and Oberman \cite{FinlayOberman2019Improved} improved their accuracy. These works treat stationary second-order equations, with no drift term and no coupled forward equation. Filtered schemes combine a base monotone scheme with a higher-order non-monotone one \cite{FroeseOberman2013,BokanowskiPicarelliReisinger2018}.

On the evacuation domain treated here, with a circular obstacle and re-entrant door jambs, a monotone grid or mesh scheme is possible but not routine, and the jambs violate the boundary smoothness that the oblique schemes above require. Particles transport the density without a mesh, but the backward equation still needs a stencil, and a scattered cloud has no mesh connectivity to build one from. A cloud that moved with the particles would also need a new stencil at every iteration.

We discretize the backward HJB equation on a fixed scattered cloud, transport the forward FP law with reflected and killed particles, and couple the two by interpolation and mass-preserving kernel density estimation (KDE). The method is a particle-in-cell architecture \cite{BirdsallLangdon1991} in which the fixed cloud takes the place of the mesh. We impose monotonicity row by row. A quadratic program (QP) returns Laplacian weights with nonnegative off-center entries, and the drift is written as a nonnegative combination of rays from the node to visible neighbors that surround it with no angular gap above a fixed cap. We call this a guarded simplicial fan. Where the physical row cannot be completed near a reflecting wall, mirror images of the neighbors across the wall are added, and next to an exit, the node takes the exit's Dirichlet value. If none of these is possible, the node is refused, which stops the run unless the node lies next to a boundary vertex, where it copies the value of a neighbor. \autoref{tab:method_comparison} places the method among the alternatives.

\begin{table}[!htbp]
    \footnotesize
    \setlength{\tabcolsep}{2.5pt}
    \centering
    \caption{Comparison of MFG numerical methods. Mass conservation: \emph{adjoint} (adjoint-consistent HJB--FP pair), \emph{transport} (characteristic-based), and \emph{particle-exact} (empirical transport, with Skorokhod reflection in the proposed method). Discrete maximum principle (DMP): \emph{upwind}, \emph{posL+fan} (positive Laplacian plus simplicial-fan upwind), and \emph{-dep.} (dependent on the chosen FEM or interpolation variant). GFDM: generalized finite difference method. $^\dagger$KDE preserves its input mass: unit mass under reflection and surviving sub-probability mass with exits.}
    \label{tab:method_comparison}
    \begin{tabular}{lcccc}
        \toprule
        Method (Ref.) & Geometry & Mass Conserv. & DMP & HJB-stab. \\
        \midrule
        FDM \cite{Achdou2010Mean,Achdou2012Mean} & Cartesian & adjoint & $\checkmark$ (upwind) & $\checkmark$ \\
        FEM \cite{Osborne2025Finite} & Simplicial & scheme-dep. & $\checkmark$ (stabilized) & $\checkmark$ \\
        Semi-Lag. \cite{Carlini2015Semi2nd} & Grid & transport & interp.-dep. & $\checkmark$ \\
        Particle \cite{DiFrancescoEtAl2022} & Meshfree & particle-exact &---&---\\
        Learning and flow matching \cite{Lauriere2022Learning,Yu2025FlowMatching} & Meshfree & not guar. &---&---\\
        \textbf{GFDM-particle} & Meshfree & particle-exact$^\dagger$ & $\checkmark$ (posL+fan) & $\checkmark$ \\
        \bottomrule
    \end{tabular}
\end{table}

\paragraph{Contributions of the current study.}
\begin{enumerate}
    \item \textbf{A monotone HJB scheme on scattered clouds.} With nonnegative off-center weights and zero row sums in each row, the Bellman residual is monotone, the matrix of each fixed policy is an M-matrix at any cell-P\'eclet number, and a discrete comparison principle holds. Upwind grids and edge-averaged or stabilized elements reach that guarantee on a mesh \cite{XuZikatanov1999,Osborne2025Finite}. Here it holds on a scattered cloud, with no separation of the cloud points.
    \item \textbf{Convergence on polygons.} On polygons, re-entrant corners included, one HJB solve at a prescribed density converges to the viscosity solution off the boundary vertices, for convex Hamiltonians of at most quadratic growth under a structure condition. A mirror row enforces normal reflection at the wall.
    \item \textbf{A coupled meshfree--particle architecture.} We build the stencil once---cloud, visibility, positive-QP rows, fan geometry and sparsity pattern---and reuse it throughout the outer iteration. Only the policy matrix is reassembled and refactorized at each Howard step, and the forward stage keeps exact Lagrangian mass accounting.
    \item \textbf{Verification against exact solutions.} Exact solutions on a floor plan test the rows of the theorem. A closed-form benchmark on a square measures rates and, against a five-point finite-difference solve, the cost of the construction, and a source-free exact MFG measures the accuracy of the fully coupled computation. On a non-convex two-exit game, we check each computed row on nine clouds.
\end{enumerate}

\section{Eulerian and Lagrangian formulations}\label{sec:dual_perspective}

Let $\Omega \subset \mbbR^d$ be a bounded Lipschitz domain, let $\mcalA \subseteq \mbbR^{d_\alpha}$ be closed, and compact unless $\mcalL/(1+\abs{b})\to\infty$ as $\abs{\alpha}\to\infty$, and take an uncontrolled constant isotropic diffusion $\sigma>0$. When the reflecting wall $\Gamma_{\mathrm{refl}}$ is nonempty, we assume that normal reflection on it defines a well-posed Skorokhod problem. The nonsmooth-domain theory for oblique reflection is developed in \cite{DupuisIshii1993,DupuisIshii2008Correction}. No convergence theorem below uses the reflected dynamics. In the interior, a representative agent with feedback control $\alpha_t = \pi(t,X_t)$ follows, for a standard $d$-dimensional Brownian motion $\ensemble{W_t}$,
\begin{equation}\label{eq:individual-dynamics}
    \diff X_t = b(X_t,m_t,\alpha_t)\diff t + \sigma\diff W_t, \qquad X_0 \sim m_0.
\end{equation}
Let $\tau_{\mathrm{abs}}:=\inf\{t\geq0:X_t\in\Gamma_{\mathrm{abs}}\}$ be the hitting time of the absorbing exit $\Gamma_{\mathrm{abs}}$, with $\tau_{\mathrm{abs}}=\infty$ when $\Gamma_{\mathrm{abs}}=\emptyset$. With reflection on $\Gamma_{\mathrm{refl}}$ and killing on $\Gamma_{\mathrm{abs}}$, the agent minimizes the stopped cost built from a running cost $\mcalL$, a terminal cost $G$, and an exit datum $g_{\mathrm{exit}}$, as follows
{\small
\begin{equation}\label{eq:stopped_cost}
    J_m(\pi)=\mbbE\bracks*{
        \int_0^{T\wedge\tau_{\mathrm{abs}}}\mcalL(X_t,m_t,\alpha_t)\diff t
        +\mathbb1_{\{\tau_{\mathrm{abs}}>T\}}G(X_T,m_T)
        +\mathbb1_{\{\tau_{\mathrm{abs}}\leq T\}}
        g_{\mathrm{exit}}(\tau_{\mathrm{abs}},X_{\tau_{\mathrm{abs}}})
    }.
\end{equation}}%
Define
\begin{equation}\label{eq:bellman_hamiltonian}
    H(x,p,m):=\sup_{\alpha\in\mcalA}
    \{-p\cdot b(x,m,\alpha)-\mcalL(x,m,\alpha)\},
\end{equation}
and fix a measurable optimal-control selector and its associated optimal drift,
\begin{align*}
    \alpha^*(x,p,m)
    &\in\argmax_{\alpha\in\mcalA}
    \{-p\cdot b(x,m,\alpha)-\mcalL(x,m,\alpha)\},\\
    b^*(x,p,m)&:=b\parens[\big]{x,m,\alpha^*(x,p,m)}.
\end{align*}

When the maximizer is unique and $H$ is differentiable in $p$, $b^*=-\nabla_pH$. The density is a probability flow under pure reflection and a sub-probability flow when particles are killed at an absorbing exit. A mean field equilibrium is a pair $(u,m)$ solving
{\small
\begin{subequations}\label{eq:MFG-system}
    \begin{empheq}[left=\empheqlbrace]{align}
        &-\partial_t u-\frac{\sigma^2}{2}\Delta u+H(x,\nabla u,m)=0, & (t,x)\in[0,T)\times\Omega, \label{eq:HJB}\\
        &\partial_t m-\frac{\sigma^2}{2}\Delta m+\diver\parens[\big]{m\,b^*(x,\nabla u,m)}=0, & (t,x)\in(0,T]\times\Omega, \label{eq:FP}\\
        &u(T,x) = G(x, m(T)), \quad m(0,x) = m_0(x). &
    \end{empheq}
\end{subequations}}%
On $\partial\Omega = \Gamma_{\mathrm{refl}} \cup \Gamma_{\mathrm{abs}}$, with $\Gamma_{\mathrm{abs}}$ relatively open in $\partial\Omega$ and outward normal $\nu$, we have
\begin{equation}\label{eq:MFG-bc}
    \begin{aligned}
        &\partial_\nu u = 0 \ \text{on}\ \Gamma_{\mathrm{refl}}, &\quad u = g_{\mathrm{exit}} \ \text{on}\ \Gamma_{\mathrm{abs}},\\
        &\parens[\Big]{\tfrac{\sigma^2}{2}\nabla m-m\,b^*(x,\nabla u,m)}\cdot\nu=0 \ \text{on}\ \Gamma_{\mathrm{refl}}, &\quad m=0 \ \text{on}\ \Gamma_{\mathrm{abs}}.
    \end{aligned}
\end{equation}

\paragraph{Solution concept.}
Viscosity solutions and boundary conditions in the viscosity sense are understood as in \cite[Secs.~7.A and~8]{CrandallIshiiLions1992}, relative to $[0,T)$, so that the inequalities include $t=0$. The exit condition holds on $[0,T)$ and the terminal condition at $t=T$, so the terminal cost and $g_{\mathrm{exit}}$ need not agree on $\Gamma_{\mathrm{abs}}$. Existence and comparison for the mixed terminal--boundary-value problem are posed off the boundary points where the outward normal is undefined or the boundary condition changes type, and off the exit at $t=T$. We assume them in general, prove comparison for a class of Hamiltonians in \autoref{app:comparison}, and obtain existence on polygons in \autoref{thm:polygon}. When they hold, $u$ denotes the unique bounded viscosity solution, and its continuity on that set is the only solution regularity the convergence theorem uses.

\paragraph{Particle representation.}
Substituting the selected optimal drift into \eqref{eq:individual-dynamics} gives the closed-loop McKean--Vlasov stochastic differential equation (SDE)
\begin{equation}\label{eq:McKean_Vlasov_SDE}
    \begin{aligned}
        \diff X_t&=b^*(X_t,\nabla u(t,X_t),m_t)\diff t+\sigma\diff W_t
        -\nu(X_t)\diff\ell_t, \qquad t<\tau_{\mathrm{abs}},\\
        m_t(B)&=\mbbP(X_t\in B,\ t<\tau_{\mathrm{abs}}),
    \end{aligned}
\end{equation}
where $B\subseteq\overline\Omega$ is any Borel set, $\ell_t$ is the reflecting-boundary local time and the process is killed at $\tau_{\mathrm{abs}}$. Under standard Lipschitz and Skorokhod well-posedness assumptions, its law is a weak solution of \eqref{eq:FP} \cite{LionsSznitman1984}. The empirical measure of $N$ particles, with a survival indicator when exits kill them, is expected to converge weakly to $m_t$, as propagation of chaos gives in the classical setting \cite{Sznitman1991}.

\section{Coupled GFDM-particle MFG solver}\label{sec:collocation}

We turn to the discrete realization of the system \eqref{eq:MFG-system}--\eqref{eq:McKean_Vlasov_SDE}. We employ a generalized finite difference method (GFDM) to approximate the spatial derivatives of the value function on a static point cloud. The collocation points $\{x_i\}_{i=1}^{N_{\mathrm{col}}}$ remain fixed throughout the computation, distinct from the moving particles $\{X_l(t)\}_{l=1}^{N_p}$ that carry probability mass for the FP equation. The two representations communicate by interpolation and KDE.

\begin{figure}[htbp]
    \centering
    \resizebox{0.90\textwidth}{!}{%
        \begin{tikzpicture}[
                node distance=2.5cm and 3.5cm,
                box/.style={rectangle, draw, rounded corners=3pt, minimum width=2.8cm, minimum height=1.2cm, align=center, font=\small, line width=0.7pt},
                eulerian/.style={box, fill=oceanicCloud3, draw=oceanicAbyss1, text=oceanicInk1},
                lagrangian/.style={box, fill=oceanicSand3, draw=oceanicCoral1, text=oceanicInk1},
                operation/.style={rectangle, draw, rounded corners=2pt, minimum width=2.2cm, minimum height=0.8cm, align=center, font=\scriptsize, fill=oceanicCloud2, draw=oceanicFoam1, line width=0.5pt, text=oceanicInk1},
                arrow/.style={->, thick, >=stealth},
                label/.style={font=\scriptsize, midway},
                oplabel/.style={font=\scriptsize, inner sep=2pt, text=oceanicSlate1},
            ]
            \node[eulerian] (hjb) {HJB Solver\\(Fixed Collocation)};
            \node[eulerian, below=1.8cm of hjb] (density) {Density $M^{(k)}$\\(Collocation Values)};
            \node[lagrangian, right=6cm of hjb] (particles) {Particle Cloud\\$\{X_l(t)\}_{l=1}^{N_p}$};
            \node[lagrangian, below=1.8cm of particles] (langevin) {SDE Integration\\(Langevin Dynamics)};
            \node[operation, right=0.7cm and 1.7cm of hjb] (interp) {$\mcalI^{\mathrm{vec}}$: Interp.};
            \node[operation, right=0.7cm and 1.7cm of density] (kde) {$\mcalR_{h_{\mathrm{KDE}}}$: KDE};
            \draw[arrow, oceanicAbyss1] (hjb) -- node[label, above, sloped] {$b_i^{*}$} (interp);
            \draw[arrow, oceanicAbyss1] (interp) -- node[label, above, sloped] {$b_h^{*}$} (particles);
            \draw[arrow, oceanicCoral1] (particles) -- node[label, right] {$\diff X_t$} (langevin);
            \draw[arrow, oceanicCoral1] (langevin) -- node[label, below, sloped] {$\{X_l(t_n)\}$} (kde);
            \draw[arrow, oceanicCoral1] (kde) -- node[label, below, sloped] {$\widetilde M^{(k+1)}$} (density);
            \draw[arrow, oceanicSlate1] (density) -- node[label, left] {$H(\cdot,M^{(k)})$} (hjb);
            \node[above=0.5cm of hjb, font=\footnotesize\bfseries, oceanicAbyss1] {Eulerian (Stable)};
            \node[above=0.5cm of particles, font=\footnotesize\bfseries, oceanicCoral1] {Lagrangian (Mass-Tracking)};
            \draw[arrow, dashed, thick, oceanicSlate1] ([xshift=-0.8cm]density.west) -- ([xshift=-0.8cm]hjb.west) node[midway, left, font=\scriptsize, text=oceanicSlate1] {fictitious play};
    \end{tikzpicture}}
    \caption{Operator structure of the GFDM-particle solver. The fixed-cloud HJB stage selects the nodal optimal control and its associated drift, interpolation passes that drift to the particles, reflected/killed SDE evolution advances the law, and KDE reconstructs the next density iterate.}
    \label{fig:architecture}
\end{figure}
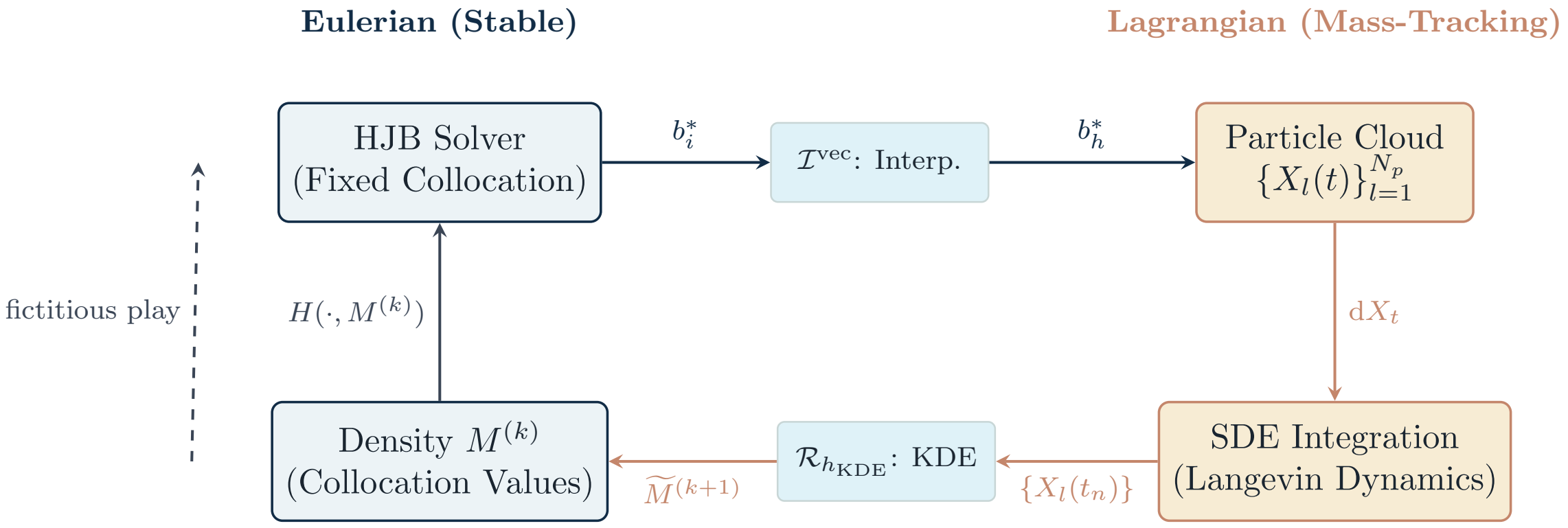

Since the FP density $m$ is the law of the McKean--Vlasov SDE \eqref{eq:McKean_Vlasov_SDE}, the FP equation is simulated by evolving particles. The coupling proceeds iteratively (\autoref{fig:architecture}): (1) solve the HJB backward on the fixed cloud for a given density iterate and select the associated nodal optimal drift, (2) evolve particles forward with its interpolation, and (3) reconstruct the next density by KDE.

\subsection{Cloud and stencil quantities}\label{sec:geometric_foundation}

We use the standard meshfree notion of fill distance \cite{Wendland2005}, localized here to the candidate neighborhoods. Write $\mcalX_h:=\{x_i\}_{i=1}^{N_{\mathrm{col}}}$ for a collocation cloud in the refinement family and set $\mcalX_h^{\mathrm{int}}:=\mcalX_h\cap\Omega$.

\begin{definition}[Stencil geometry]\label{def:stencil_geometry}
    Let $x_i\in\mcalX_h$, and let its candidate neighborhood satisfy
    $x_i\in\mcalV_i\subset\mcalX_h\cap B(x_i,\delta_i)\cap\overline\Omega$, with $(x_i,x_j)\subset\Omega$ for each $x_j\in\mcalV_i\setminus\{x_i\}$.
    The local fill distance, taken over the points of the ball that $x_i$ sees, is
    \begin{equation*}
        h_i:=\sup_{\substack{x\in B(x_i,\delta_i)\cap\Omega\\(x_i,x)\subset\Omega}}
        \min_{x_j\in\mcalV_i}\abs{x-x_j}.
    \end{equation*}
\end{definition}

\begin{definition}[Visibility]\label{def:visibility}
    An edge $[x_i,x_j]$ with $x_i,x_j\in\overline\Omega$ is \emph{visible} if the open segment $(x_i,x_j)$ lies in $\Omega$, and a candidate neighborhood $\mcalV_i$ is \emph{visibility-safe} if all edges from $x_i$ to members of $\mcalV_i$ are visible. Write $\mcalV_i^{\mathrm{vis}}(\delta):=\{x_j\in\mcalX_h\cap B(x_i,\delta)\cap\overline\Omega:[x_i,x_j]\text{ visible}\}$ for the largest such neighborhood of radius $\delta$.
\end{definition}

Visibility filtering discards the remaining edges at stencil assembly, so no stencil couples nodes across solid material. For a refinement family, define $h:=\max_{x_i\in\mcalX_h^{\mathrm{int}}}h_i$. Since only points visible from $x_i$ enter $h_i$, a region hidden behind an obstacle inflates neither $h_i$ nor $h$.

\begin{definition}[Collocation family with local support]\label{def:local_support}
    A refinement family has \emph{local support} if there exist $c_\delta,C_\delta>0$, independent of $h$, such that
    \begin{equation*}
        c_\delta h\leq\delta_i\leq C_\delta h,
        \qquad x_i\in\mcalX_h^{\mathrm{int}},
    \end{equation*}
    and $\sup_{x\in\Omega}\operatorname{dist}(x,\mcalX_h)\to0$ as $h\to0$.
\end{definition}

We ask for no separation of the nodes, since the two-dimensional Sobol sequence is not quasi-uniform \cite{Goda2024Sobol}. The estimates below rest on the nonnegativity and the exact moments of the weights, on the aperture cap and on the locality of the stencils.

\subsection{Row construction}\label{sec:row_construction}

\paragraph{The positive-Laplacian program.}

The diffusion operator is the sparse nodal map
\begin{equation}\label{eq:gfdm_weights}
    (\Delta_hu)_i:=\sum_{j\in\mcalV_i}L_{ij}u_j.
\end{equation}
Write $\msfP_k$ for the space of real polynomials of degree at most $k$. Seibold's decoupled monotone Laplacian construction \cite{Seibold2008Minimal}, extended to non-divergence-form operators by Ye and Tian \cite{YeTian2023}, imposes nonnegativity on the weights of a row that is exact on $\msfP_2$,
\begin{equation}\label{eq:qp_problem}
    \begin{aligned}
        \min_{L_{i,\cdot}} \quad & \norm{L_{i,\cdot}}_2^2 \\
        \text{s.t.} \quad & \sum_jL_{ij}q(x_j)=\Delta q(x_i)\ \text{ for all }q\in\msfP_2,
        \quad L_{ij} \geq 0 \ (j \neq i),
    \end{aligned}
\end{equation}
optimizing directly over the physical stencil row $L_{i,\cdot}:=(L_{ij})_j$ (center and neighbors). Taking $q\equiv1$ gives $\sum_jL_{ij}=0$. The quadratic objective of the QP \eqref{eq:qp_problem} selects a minimal-norm feasible row.

\paragraph{The simplicial-fan advection operator.}
The interior and reflected fan geometry is illustrated in \autoref{fig:fan_schematic}.
Throughout, $d=2$. For a local density value $m_i$ and a fixed control $\alpha\in\mcalA$, we discretize the advective term $b(x_i,m_i,\alpha)\cdot\nabla u(x_i)$ by a \emph{simplicial fan}. In the obstacle-free interior, ordered stencil rays $z_{ij}:=x_j-x_i$ generate a complete face-to-face complex of two-ray cones covering $\mbbR^2$. Visibility filtering may leave the physical candidate fan, built from visible cloud points alone, open near an obstacle. Where a reflecting face is within reach, reflected rays complete it (\autoref{def:row_classification}). The retained physical or reflected fan is geometry-only and fixed before any policy solve. On a complete fan, the drift of a control lies in a retained sector. Let $\lambda_{ij}^{\alpha}$ denote its nonnegative conic coordinates on the two bracketing rays, extended by zero to all other rays, which depend on $m_i$ only through $b$, and define
\begin{equation}\label{eq:fan_advection}
    \begin{aligned}
        &b(x_i,m_i,\alpha)=\sum_j\lambda_{ij}^{\alpha}z_{ij},
        &\quad\lambda_{ij}^{\alpha}\geq0, \\
        &\sum_j\lambda_{ij}^{\alpha}\abs{z_{ij}}\leq C_\theta\abs{b(x_i,m_i,\alpha)},
        &\quad(A_{\mathrm{fan}}^\alpha u)_i
        :=\sum_j\lambda_{ij}^{\alpha}(u_j-u_i).
    \end{aligned}
\end{equation}
The existence of these coordinates and the stated bound follow from \autoref{lem:fan_existence}. Because the representation is exact on affine data, $A_{\mathrm{fan}}^\alpha$ is first-order consistent with $b(x_i,m_i,\alpha)\cdot\nabla$. Because its weights are nonnegative, $-A_{\mathrm{fan}}^\alpha$ contributes non-positive off-diagonals to the implicit matrix. A single operator is therefore both monotone and first-order consistent. In two dimensions, for a fixed drift direction, the row \eqref{eq:fan_advection} is the monotone first-derivative scheme that Finlay and Oberman \cite[Def.~1 and Lem.~1]{FinlayOberman2019Improved} write in barycentric coordinates on a point cloud, which is degenerate elliptic in the sense of Oberman \cite{Oberman2006Convergent}. They use it as a building block for second derivatives, in equations without a drift term. Here it discretizes the drift of the Bellman operator and enters one residual with the positive Laplacian.

\begin{figure}[htbp]
    \centering
    \begin{tikzpicture}[scale=1.12, >={Latex[length=1.4mm,width=1.0mm]}]
        \begin{scope}
            \draw[dashed, oceanicSlate1, opacity=0.45] (0,0) circle (1.5);
            \foreach \aa/\bb in {18/64,64/137,137/191,191/267,267/316,316/378} {
                \filldraw[fill=oceanicCloud3, fill opacity=0.25, draw=white, line width=0.6pt]
                (0,0) -- (\aa:0.88) arc (\aa:\bb:0.88) -- cycle;
            }
            \fill[oceanicCoral1, opacity=0.18] (0,0) -- (18:0.88) arc (18:64:0.88) -- cycle;
            \foreach \a/\r in {18/1.43, 64/1.03, 137/1.48, 191/0.92, 267/1.38, 316/1.12} {
                \draw[thin, oceanicSlate1, opacity=0.7] (0,0) -- (\a:\r);
                \fill[oceanicAbyss1] (\a:\r) circle (1.4pt);
            }
            \draw[thick, oceanicCoral1] (0,0) -- (18:1.43);
            \draw[thick, oceanicCoral1] (0,0) -- (64:1.03);
            \node[oceanicAbyss1, font=\scriptsize] at (18:1.67) {$z_{i,j}$};
            \node[oceanicAbyss1, font=\scriptsize] at (64:1.39) {$z_{i,j+1}$};
            \draw[->, very thick, oceanicCoral1] (0,0) -- (36:0.92)
            node[right = -1pt, font=\small, oceanicCoral1] {$b$};
            \draw[thin, oceanicCoral1] (55:0.68) -- (1.15,1.05) -- (1.28,1.05);
            \node[oceanicCoral1, anchor=west, font=\scriptsize] at (1.31,1.05) {$C_j$};
            \fill[oceanicAbyss1] (0,0) circle (1.7pt) node[below left=-1pt, font=\small] {$x_i$};
            \node[oceanicAbyss1, font=\footnotesize] at (0,-2.00)
            {$\bigcup_{j=1}^{6}C_j=\mbbR^2$};
            \node[oceanicSlate1, font=\small] at (0,-2.45) {(a) Interior: complete fan};
        \end{scope}
        \begin{scope}[xshift=5.6cm]
            \coordinate (pOne) at (0.82,-1.15);
            \coordinate (pTwo) at (0.78,-0.25);
            \coordinate (pThree) at (1.18,0.63);
            \coordinate (pFour) at (0.15,1.05);
            \coordinate (gOne) at (-1.52,-1.15);
            \coordinate (gTwo) at (-1.48,-0.25);
            \coordinate (gThree) at (-1.88,0.63);
            \coordinate (gFour) at (-0.85,1.05);
            \fill[oceanicSlate1, opacity=0.22] (-2.05,-1.7) rectangle (-0.35,1.7);
            \draw[oceanicSlate1, thick] (-0.35,-1.7) -- (-0.35,1.7);
            \node[oceanicSlate1, font=\scriptsize] at (-1.48,1.48) {obstacle};
            \foreach \p/\g in {pOne/gOne, pTwo/gTwo, pThree/gThree, pFour/gFour} {
                \draw[dotted, oceanicPlum3, opacity=0.55] (\p) -- (\g);
            }
            \draw[|<->|, thin, oceanicPlum3, opacity=0.62] (-0.85,1.24) -- (-0.35,1.24)
            node[midway, above=-2pt, font=\scriptsize] {$s_{ij}$};
            \draw[|<->|, thin, oceanicPlum3, opacity=0.62] (-0.35,1.24) -- (0.15,1.24)
            node[midway, above=-2pt, font=\scriptsize] {$s_{ij}$};
            \fill[oceanicSlate1, opacity=0.10] (0,0) -- (82:0.9) arc (82:305:0.9) -- cycle;
            \draw[thin, oceanicSlate1, opacity=0.70]
            (-0.73,-1.25) -- (-0.47,-0.84) -- (250:0.58);
            \node[oceanicSlate1, font=\scriptsize, align=center] at (-0.95,-1.43) {occluded\\directions};
            \filldraw[fill=oceanicCloud3, fill opacity=0.42, draw=white, line width=0.8pt]
            (0,0) -- (-55:0.78) arc (-55:-18:0.78) -- cycle;
            \filldraw[fill=oceanicSand3, fill opacity=0.35, draw=white, line width=0.8pt]
            (0,0) -- (-18:0.78) arc (-18:28:0.78) -- cycle;
            \filldraw[fill=oceanicCloud2, fill opacity=0.32, draw=white, line width=0.8pt]
            (0,0) -- (28:0.78) arc (28:82:0.78) -- cycle;
            \fill[oceanicCoral1, opacity=0.18] (0,0) -- (-18:0.78) arc (-18:28:0.78) -- cycle;
            \foreach \g in {gOne,gTwo,gThree,gFour} {
                \draw[densely dashed, oceanicMoss2, opacity=0.55] (0,0) -- (\g);
                \fill[oceanicMoss2, opacity=0.85] (\g) circle (1.2pt);
            }
            \draw[dashed, oceanicPlum3, opacity=0.62] (pThree) -- (gThree);
            \draw[thin, oceanicPlum3, opacity=0.62] (-0.35,0.75) -- (-0.23,0.75) -- (-0.23,0.63);
            \foreach \p in {pOne,pTwo,pThree,pFour} {
                \draw[thin, oceanicSlate1, opacity=0.7] (0,0) -- (\p);
                \fill[oceanicAbyss1] (\p) circle (1.4pt);
            }
            \draw[->, very thick, oceanicCoral1] (0,0) -- (8:0.9) node[right=-1pt, font=\small, oceanicCoral1] {$b$};
            \draw[thin, oceanicCoral1] (20:0.72) -- (1.30,0.27) -- (1.43,0.27);
            \node[oceanicCoral1, anchor=west, font=\scriptsize] at (1.46,0.27) {$C_2$};
            \node[oceanicAbyss1, anchor=south west, font=\scriptsize] at (1.20,0.67) {$x_j$};
            \node[oceanicMoss2, font=\scriptsize] at (-1.78,0.88) {$x_j^{\mathrm{ghost}}$};
            \fill[oceanicAbyss1] (0,0) circle (1.7pt) node[below left=-1pt, font=\small] {$x_i$};
            \node[oceanicAbyss1, font=\footnotesize] at (0.48,-2.00)
            {$\bigcup_{j=1}^{3}C_j\ne\mbbR^2$};
            \node[oceanicSlate1, font=\small] at (-0.3,-2.45) {(b) Physical fan open: reflect to complete};
        \end{scope}
    \end{tikzpicture}
    \caption{Guarded-fan geometry. Adjacent retained rays define cones $C_j$, and coral marks the cone containing $b$. \textbf{(a)} A complete interior fan tiles $\mbbR^2$. \textbf{(b)} Visibility leaves the physical fan open. Each ghost is a physical source reflected in the face, at equal distance $s_{ij}$ on either side (plum, dotted). The green rays to those ghosts complete the fan, and their weights fold back onto the physical sources.}
    \label{fig:fan_schematic}
\end{figure}

Let $\mscrS_i$ denote the adjacent ray pairs retained by the guarded fan at $x_i$.

\begin{definition}[Complete fans]\label{def:fan_complete}
    Fix an aperture cap $\theta_+\in(0,\pi)$. The fan at $x_i$ is \emph{complete} if the cones $\operatorname{cone}(z_{i,j},z_{i,j+1})$, $(j,j+1)\in\mscrS_i$, cover $\mbbR^2$ and $0<\angle(z_{i,j},z_{i,j+1})\leq\theta_+$ for each retained adjacent pair.
\end{definition}

\begin{lemma}[Simplicial-fan representation]\label{lem:fan_existence}
    Let $d=2$ and suppose the fan at $x_i$ is complete. Then every $v\in\mbbR^2$ admits coefficients satisfying
    \begin{equation*}
        v=\sum_j\lambda_{ij}z_{ij},
        \qquad
        \lambda_{ij}\geq0,
        \qquad
        \sum_j\lambda_{ij}\abs{z_{ij}}\leq C_\theta\abs{v},
    \end{equation*}
    where $C_\theta:=1/\cos(\theta_+/2)$.
\end{lemma}

\begin{proof}[Sketch]
    Every nonzero $v$ lies in one adjacent two-ray cone of a complete fan, and on a cone of aperture $\theta\leq\theta_+$ the law of cosines bounds the length-weighted coefficients by $\abs{v}/\cos(\theta/2)$. See \autoref{app:proof_fan_existence} for the details.
\end{proof}

The constructed operator installs only complete physical or reflected fans on dynamic rows, so every control is represented on each assembled row.

\paragraph{Positive-row construction with reflecting closure.}

A node whose visible pool lies in an open half-plane admits no nonnegative Laplacian row over physical points \cite[Thm.~8]{Seibold2008Minimal}. Near a reflecting face, ghost images restore the missing side. For a mirror row, let $\mscrJ_i^{\mathrm{ghost}}$ denote the physical sources that receive reflected ghosts, write $x_j^{\mathrm{ghost}}$ for the ghost of $x_j$, and let $s_{ij}\geq0$ be the inward distance of $x_j$, $j\in\mscrJ_i^{\mathrm{ghost}}$, from its selected planar face. Denote the augmented-row weight at $x_j^{\mathrm{ghost}}$ by $\widetilde L^{\mathrm{ghost}}_{ij}$, extended by zero for $j\notin\mscrJ_i^{\mathrm{ghost}}$. In addition to the positive $\msfP_2$ moments, the QP \eqref{eq:qp_problem} imposes the linear discrete-obliqueness constraint
\begin{equation}\label{eq:bc_gamma_main}
    \widehat\gamma_{i,h}:=2\delta_i\sum_{j\in\mscrJ_i^{\mathrm{ghost}}}
    \widetilde L^{\mathrm{ghost}}_{ij}s_{ij}\geq\gamma_0>0,
\end{equation}
where $\delta_i$ is the mirror-pool radius and $\gamma_0$ is independent of $h$. The mirror row is the folded row $L_{ij}:=\widetilde L_{ij}+\widetilde L^{\mathrm{ghost}}_{ij}$, with $\widetilde L_{ij}$ the augmented-row weight at $x_j$, and the reflected guarded fan folds each ghost ray onto its physical source in the same way.

\begin{definition}[Row classification]\label{def:row_classification}
    Fix a reach factor $\kappa>0$. At each interior node $x_i\in\mcalX_h^{\mathrm{int}}$, the construction first builds the physical pool $\mcalV_i^{\mathrm{vis}}(\delta_i)$ of \autoref{def:visibility} and its guarded fan, and then returns the first of the following that applies:
    \begin{enumerate}
        \item[(P)] a \emph{partial differential equation (PDE) row}, if the physical fan is complete and \eqref{eq:qp_problem} is feasible on the initial or visibility-safe enlarged physical stencil within the pool.
        \item[(D)] an \emph{absorbing-Dirichlet identity row}, if the boundary face nearest to $x_i$ is an absorbing exit at distance at most $\kappa\delta_i$.
        \item[(R)] a \emph{mirror row}, if a reflecting face lies within distance $\kappa\delta_i$ of $x_i$, the folded fan over the physical rays and the rays reflected across such faces is complete, and the augmented weights satisfy \eqref{eq:bc_gamma_main}.
        \item[(X)] a \emph{refusal}: the construction is undefined at that node.
    \end{enumerate}
    We sample no node on $\Gamma_{\mathrm{refl}}$. Nodes sampled on $\Gamma_{\mathrm{abs}}$ are not classified and take absorbing-Dirichlet identity rows with the exit datum.
\end{definition}

\begin{algorithm}[htbp]
    \footnotesize
    \SetAlgoLined
    \DontPrintSemicolon
    \KwIn{$\mcalX_h$, pool radii $\delta_i$, initial stencils $\mcalV_i^0\subseteq\mcalV_i^{\mathrm{vis}}(\delta_i)$, search radii $r_1<r_2<\cdots$, $\kappa>0$, $\theta_+\in(0,\pi)$, $\gamma_0>0$, the set $\Gamma_{\mathrm{c}}$ of boundary vertices and junctions, $C_{\mathrm{v}}>0$.}
    \KwOut{The type of every $x_i\in\mcalX_h^{\mathrm{int}}$ in \autoref{def:row_classification}, and its row.}

    \BlankLine
    \ForEach{$x_i\in\mcalX_h^{\mathrm{int}}$}{
        $\mcalV\leftarrow\mcalV_i^{\mathrm{vis}}(\delta_i)$, \quad $Z\leftarrow$ a complete fan over $\mcalV\cap B(x_i,r_k)$ with $k$ least, or $\emptyset$\;
        \lIf{$Z\neq\emptyset$ and \eqref{eq:qp_problem} is feasible on $\mcalV_i^0$ or on $\mcalV$}{type (P)}
        \lElseIf{$\operatorname{dist}(x_i,\Gamma_{\mathrm{abs}})=\operatorname{dist}(x_i,\partial\Omega)\leq\kappa\delta_i$}{type (D)}
        \uElseIf{$\operatorname{dist}(x_i,\Gamma_{\mathrm{refl}})\leq\kappa\delta_i$}{
            $\mcalV^{\mathrm{g}}\leftarrow\mcalV\cup\{x_j^{\mathrm{ghost}}:x_j\in\mcalV\}$, \quad $Z\leftarrow$ a complete fan over $\mcalV^{\mathrm{g}}$, or $\emptyset$\;
            \lIf{$Z\neq\emptyset$ and \eqref{eq:qp_problem} with \eqref{eq:bc_gamma_main} is feasible on $\mcalV^{\mathrm{g}}$}{type (R)}
            \lElse{type (X)}
        }
        \lElse{type (X)}
    }
    \ForEach{$x_i$ of type (X)}{
        \lIf{$\operatorname{dist}(x_i,\Gamma_{\mathrm{c}})\leq C_{\mathrm{v}}h$}{copy row \eqref{eq:copy_row}}
        \lElse{stop}
    }
    \caption{Row construction.}
    \label{alg:row_construction}
\end{algorithm}

The type of every row is fixed before the first policy solve, by \autoref{alg:row_construction}. A node with no admissible row stops the run rather than falling back to a non-monotone row, unless it lies within a fixed multiple of the fill distance of a boundary vertex, where it copies the value of a nearby dynamic row.

\subsection{Discrete HJB system and M-matrix structure}\label{sec:monotonicity_QP_GFDM_HJB}

\paragraph{HJB time discretization and discrete Hamiltonian.}

Let $\Delta t:=T/N_t$ and $t_n:=n\Delta t$, and let $\pi_h$ denote nodal sampling on the collocation cloud $\mcalX_h$: for a scalar spatial function $v$,
\begin{equation*}
    \pi_hv:=\parens[\big]{v(x_j)}_{j=1}^{N_{\mathrm{col}}}.
\end{equation*}
For a time-dependent field, we apply $\pi_h$ at each $t_n$.
For a smooth test function and a prescribed density field, the continuous HJB residual is
\begin{equation*}
    F[\varphi;m]
    :=-\partial_t\varphi-\frac{\sigma^2}{2}\Delta\varphi
    +H(x,\nabla\varphi,m).
\end{equation*}
The density enters through its local value, so we also write $F[\varphi;\mu](t,x)$ with a number $\mu$ in place of $m(t,x)$. For a nodal vector $W\in\mbbR^{N_{\mathrm{col}}}$, with one value per collocation point, define the spatial discrete Hamiltonian at $x_i$ by
\begin{equation}\label{eq:discrete_hamiltonian}
    H_{i,h}[W;m_i]
    :=\sup_{\alpha\in\mcalA}
    \ensemble*{-(A_{\mathrm{fan}}^\alpha W)_i-\mcalL(x_i,m_i,\alpha)}.
\end{equation}
When this supremum is attained, select the discrete optimal control as
\begin{equation}\label{eq:discrete_optimal_control}
    \alpha_{i,h}^{*}[W;m_i]
    \in\argmax_{\alpha\in\mcalA}
    \ensemble*{-(A_{\mathrm{fan}}^\alpha W)_i-\mcalL(x_i,m_i,\alpha)}.
\end{equation}
The associated discrete optimal drift is
\begin{equation*}
    b_{i,h}^{*}[W;m_i]
    :=b\parens[\big]{x_i,m_i,\alpha_{i,h}^{*}[W;m_i]}.
\end{equation*}
We suppress the arguments $h$, $W$, $m_i$ in $b_i^{*}$ whenever they are fixed and unambiguous. Thus, $\alpha$ denotes a generic control, $\alpha_{i,h}^{*}$ the selected discrete optimal control, and $b_{i,h}^{*}$ its associated optimal drift. The last two coincide only in the identity-drift case $b(x,m,\alpha)=\alpha$. If $x_i$ has a complete physical fan, then for all $\varphi\in\msfP_1$ and all density values $m_i$ for which the continuous Hamiltonian is finite, we have
\begin{equation*}
    H_{i,h}[\pi_h\varphi;m_i]
    =H(x_i,\nabla\varphi,m_i).
\end{equation*}

A \emph{discrete density history} is a family $M=(M^n)_{n=0}^{N_t}\in(\mbbR^{N_{\mathrm{col}}})^{N_t+1}$, where $M_i^n$ denotes the density value at the node $x_i$ and the time $t_n$. It is the input of one HJB solve, either the density iterate $M^{(k)}$ at outer iteration $k$ or the nodal sample of a prescribed density field. For a nodal history $W=(W^n)_{n=0}^{N_t}\in(\mbbR^{N_{\mathrm{col}}})^{N_t+1}$ and a discrete density history $M$, define the implicit Euler residual by
\begin{equation}\label{eq:F_fan}
    F_{h,\Delta t}^{n,i}(W;M)
    :=-\frac{W_i^n-W_i^{n-1}}{\Delta t}
    -\frac{\sigma^2}{2}(\Delta_hW^{n-1})_i
    +H_{i,h}
    [W^{n-1};M_i^n].
\end{equation}
A residual labeled by a space--time point has the time level of the value it determines. Thus the residual $F_{h,\Delta t}(W;M)(t_{n-1},x_i)$ is $F_{h,\Delta t}^{n,i}(W;M)$.
Let $\mscrI_h^{\mathrm{dyn}}$ be the space--time indices of the interior and mirror rows and $\mscrI_h^{\mathrm{data}}$ the terminal and absorbing-Dirichlet indices. A node next to a boundary vertex that admits no row has neither a residual nor an imposed datum, and we collect these indices in a third class $\mscrI_h^{\mathrm{cl}}$, each equipped with the \emph{copy row}
\begin{equation}\label{eq:copy_row}
    U_i^n-U_{\iota(i)}^n=0,
    \qquad(n,i)\in\mscrI_h^{\mathrm{cl}},
\end{equation}
where $\iota(i)$ indexes a dynamic row within a fixed multiple of the fill distance, $\abs{x_{\iota(i)}-x_i}\leq C_{\mathrm{cl}}h$. A copy row annihilates constants and has the off-diagonal sign of the other two classes. Dynamic rows may use copied nodes as stencil sources, that is, $L_{ji}>0$ or $\lambda_{ji}^\alpha>0$ is allowed for a dynamic row $j$ and a copied node $i$. When no node is copied, $\mscrI_h^{\mathrm{cl}}=\emptyset$ and all statements below reduce to the two-class case. Given $M^{(k)}$, let $g_h^{(k)}$ denote the imposed data, with
\begin{equation*}
    g_{h,i}^{(k),n}
    :=
    \begin{cases}
        G(x_i,M^{(k),N_t}),&n=N_t,\\
        g_{\mathrm{exit}}(t_n,x_i),&n<N_t,\ x_i\text{ carrying an absorbing-Dirichlet row}.
    \end{cases}
\end{equation*}
Such a row is installed at a node near an absorbing exit, and the node need not lie on $\Gamma_{\mathrm{abs}}$. The exit datum is extended to that boundary layer. The discrete HJB system is
\begin{equation}\label{eq:discrete_hjb_system}
    \begin{aligned}
        F_{h,\Delta t}^{n,i}(U^{(k+1)};M^{(k)})&=0,
        &&(n,i)\in\mscrI_h^{\mathrm{dyn}},\\
        U_i^{(k+1),n}&=g_{h,i}^{(k),n},
        &&(n,i)\in\mscrI_h^{\mathrm{data}},\\
        U_i^{(k+1),n}-U_{\iota(i)}^{(k+1),n}&=0,
        &&(n,i)\in\mscrI_h^{\mathrm{cl}}.
    \end{aligned}
\end{equation}
We hold the prescribed density iterate fixed through the backward sweep. Here $M_i^{(k),n}$ is the local density value sampled at the right endpoint of $[t_{n-1},t_n]$, a first-order time quadrature consistent with implicit Euler. The residual \eqref{eq:F_fan} is used by both the algorithm and the convergence proof.

With $b$ and $\mcalL$ continuous in $\alpha$, the standing assumption of a compact control set, or of a running cost superlinear in the drift, makes the discrete Bellman supremum attained. The comparison principle itself needs only $\varepsilon$-maximizers.

For a policy $\alpha=(\alpha_i)_i$, the fully implicit matrix acts by
\begin{equation}\label{eq:implicit_matrix}
    (A^\alpha W)_i:=
    \begin{cases}
        \dfrac{W_i}{\Delta t}-\dfrac{\sigma^2}{2}(\Delta_hW)_i-(A_{\mathrm{fan}}^\alpha W)_i, & x_i\text{ of type (P) or (R)},\\
        W_i, & x_i\text{ absorbing-Dirichlet},\\
        W_i-W_{\iota(i)}, & x_i\text{ copied},
    \end{cases}
\end{equation}
and a fixed-policy step of \eqref{eq:discrete_hjb_system} is a linear system with this matrix.

``M-matrix'' means nonsingular M-matrix, with the characterizations of \cite[Chap.~6]{BermanPlemmons1994}.

\begin{theorem}[Algebraic monotonicity of the constructed operator]\label{thm:constructed_operator_monotone}
    Suppose every dynamic row is of type (P) or (R) in the sense of \autoref{def:row_classification}, with every ghost image, across one or several faces, being folded onto its physical source. Then:
    \begin{enumerate}
        \item[(i)] every non-Dirichlet row constructed in this way satisfies $L_{ij} \geq 0$ $(j \neq i)$ and $\sum_j L_{ij} = 0$, and hence the spatial generator $\Delta_h + A_{\mathrm{fan}}^\alpha$ is Metzler with zero row sum on each such row, and
        \item[(ii)] for any fixed policy $\alpha$, the matrix $A^\alpha$ of \eqref{eq:implicit_matrix} is a nonsingular M-matrix for all $\Delta t>0$.
    \end{enumerate}
\end{theorem}

\begin{proof}[Sketch]
    Feasibility of \eqref{eq:qp_problem} gives nonnegative off-center coefficients and zero row sum. The augmented mirror program has the same properties, which ghost folding preserves. The fan adds another Metzler, zero-row-sum contribution. For a fixed policy, the dynamic rows of $A^\alpha$ are strictly diagonally dominant $Z$-matrix rows, while the data rows are identities. Eliminating the copied unknowns leaves a strictly diagonally dominant $Z$-matrix, so the whole matrix is a nonsingular M-matrix. The folding and block details are in \autoref{app:proof_constructed_operator}.
\end{proof}

\paragraph{HJB nonlinear inner solve: Howard's policy iteration.}

We solve each backward step of \eqref{eq:discrete_hjb_system} by exact-policy Howard iteration \cite{BokanowskiMarosoZidani2009}. All the benchmarks below use the quadratic drift-control instance $b(x,m,\alpha)=\alpha$, $\mcalL=\frac{c}{2}\abs{\alpha}^2+V(t,x)+f(m)$ with constant $c>0$, which is coercive, so that $\mcalA=\mbbR^2$. Then $H_{i,h}$ is piecewise quadratic over the retained sectors, and policy improvement therefore reduces to finitely many Euclidean projections and each policy-evaluation solve uses an M-matrix. This closed-form policy update is specific to the instance, while other Bellman data require an exact or convergent policy optimizer. We collect the update equations and the descent statement in \autoref{app:howard_solver}.

\subsection{Coupling with the particles}\label{sec:collocation_coupling}

\paragraph{Reflected-SDE particle scheme.}

The particle step discretizes the reflected-and-killed dynamics \eqref{eq:McKean_Vlasov_SDE}. With the interpolated discrete optimal drift, Euler--Maruyama \cite{Kloeden1992} gives
\begin{equation}\label{eq:euler_maruyama}
    \begin{aligned}
        X_l(t_{n+1}) &= R_\Gamma\parens[\big]{X_l(t_n),\,Y_n^l},\\
        Y_n^l &:= X_l(t_n) + b^{*}_h(X_l(t_n)) \Delta t + \sigma \sqrt{\Delta t} Z_n^l,
        \qquad Z_n^l \sim \mcalN(0, \Id),
    \end{aligned}
\end{equation}
where $R_\Gamma(x,y)$ resolves the first contact of the segment $[x,y]$ with $\Gamma_{\mathrm{refl}}$ in closed form and reflects the remaining displacement about the contacted normal, iterating within the step. A normal projection onto $\overline\Omega$ is retained as a safeguard and is inactive once the contact is resolved. A trajectory is killed when its proposed segment first meets $\Gamma_{\mathrm{abs}}$. On the straight exits this is a test at the end of each step, so an excursion across $\Gamma_{\mathrm{abs}}$ within a step goes undetected. On smooth domains, killing tested only at the steps has a weak error of order $\Delta t^{1/2}$ \cite[Thm.~2.4]{Gobet2000}.

\paragraph{Cloud--particle transfer.}

At each time slice, let $\mu_t^{N_p}:=N_p^{-1}\sum_{\tau_l>t}\delta_{X_l(t)}$ denote the surviving empirical sub-probability measure. Write $\mcalR_{h_{\mathrm{KDE}}}$ for a reconstruction that maps it to a nodal density with the same total mass. On geometries admitting an exact image sum, it is the domain-kernel operator
\begin{equation}\label{eq:kde_operator}
    \mcalK_{h_{\mathrm{KDE}}}[\mu](x)
    :=\int_{\overline\Omega}K^\Omega_{h_{\mathrm{KDE}}}(x,y)\diff\mu(y),
    \qquad
    \int_{\overline\Omega}K^\Omega_{h_{\mathrm{KDE}}}(x,y)\diff x=1
    \quad\text{for all }y\in\overline\Omega.
\end{equation}
Consequently, $\int_{\overline\Omega}\mcalK_{h_{\mathrm{KDE}}}[\mu]=\mu(\overline\Omega)$. On obstacle geometry, we instead use the ghost-image reconstruction $\mcalK^{\mathrm{gh}}_{h_{\mathrm{KDE}}}$ of \autoref{app:kde_construction}, which truncates the image sum and rescales the result to the input mass.

The reverse transfer applies the partition-of-unity interpolant $\mcalI^{\mathrm{vec}}[\mbfv](x):=\sum_i\psi_i(x)v_i$, with $\sum_i\psi_i(x)=1$, to the nodal optimal drift.

\autoref{alg:Hybrid-Eulerian-Lagrangian} combines these stages into the coupled solver. It builds the rows once, and at each outer iteration it solves the HJB system backward at the current density, moves the particles forward under the resulting drift, reconstructs their density by kernel density estimation, and averages it into the next density iterate.

\begin{algorithm}[htbp]
    \footnotesize
    \SetAlgoLined
    \DontPrintSemicolon
    \KwIn{Initial density $m_0$, terminal cost $G$, time horizon $T$, parameters $(N_t,N_p,N_{\mathrm{col}},\sigma)$, number of outer iterations $K$, damping schedule $(\eta_k)_{k\geq0}\subset(0,1]$.}
    \KwOut{Final iterate $(U^{(K)},M^{(K)})$.}

    \BlankLine
    Build the rows by \autoref{alg:row_construction}\;
    Sample particles from $m_0$ and initialize the density history\;

    \BlankLine
    \For{$k\leftarrow0$ \KwTo $K-1$}{
        Solve \eqref{eq:discrete_hjb_system} backward by exact-policy Howard iteration\;
        Evolve particles by \eqref{eq:euler_maruyama}\;
        Reconstruct $\widetilde M^{(k+1)}$ by $\mcalR_{h_{\mathrm{KDE}}}$ and update
        $M^{(k+1)}\leftarrow(1-\eta_k)M^{(k)}+\eta_k\widetilde M^{(k+1)}$\;
        Record the outer fixed-point residual $\norm{\widetilde M^{(k+1)}-M^{(k)}}_1$\;
    }
    \caption{GFDM-particle MFG solver.}
    \label{alg:Hybrid-Eulerian-Lagrangian}
\end{algorithm}

With $\eta_k=1/(k+1)$, \autoref{alg:Hybrid-Eulerian-Lagrangian} is a stochastic particle realization of MFG fictitious play \cite{Cardaliaguet2017}, and $M^{(k)}$ is the running average of the past responses.

\section{Comparison and convergence at a prescribed density}\label{sec:hjb_convergence}

In this section, we prove a discrete comparison principle and stability, which hold at each density iterate, and then show the convergence of one HJB solve at a prescribed density, in the interior and up to planar reflecting segments and exit windows.

For the comparison theory, we fix a discrete density history $M=(M^n)_{n=0}^{N_t}$ as in \eqref{eq:F_fan}. Denote the corresponding imposed data by $g_h=(g_{h,i}^n)$, indexed by $(n,i)\in\mscrI_h^{\mathrm{data}}$.

\begin{definition}[Discrete subsolution and supersolution]\label{def:discrete_subsolution}
    A grid function $\underline U$ is a discrete subsolution if
    \begin{align*}
        F_{h,\Delta t}^{n,i}(\underline U;M)&\leq0,
        \quad (n,i)\in\mscrI_h^{\mathrm{dyn}},\\
        \underline U_i^n&\leq g_{h,i}^n,
        \quad (n,i)\in\mscrI_h^{\mathrm{data}},\\
        \underline U_i^n-\underline U_{\iota(i)}^n&\leq0,
        \quad (n,i)\in\mscrI_h^{\mathrm{cl}}.
    \end{align*}
    A grid function $\overline U$ is a discrete supersolution if all three inequalities above, with $\underline U$ replaced by $\overline U$, are reversed.
\end{definition}

\begin{theorem}[P\'eclet-free discrete comparison principle]\label{thm:upwind_comparison}
    Let $H$ have the Bellman form \eqref{eq:bellman_hamiltonian}, with finite suprema in \eqref{eq:discrete_hamiltonian}. Suppose every dynamic row has a complete fan and a Laplacian row with $L_{ij}\geq0$ for $j\neq i$ and $\sum_jL_{ij}=0$, absorbing-Dirichlet nodes are imposed by identity rows, and any copied nodes are assigned the copy rows \eqref{eq:copy_row}.
    If $\underline U$ is a discrete subsolution and $\overline U$ is a discrete supersolution in the sense of \autoref{def:discrete_subsolution}, then
    \begin{equation*}
        \underline U_i^n\leq\overline U_i^n,
        \qquad x_i\in\mcalX_h,\quad n=0,\ldots,N_t.
    \end{equation*}
\end{theorem}

\begin{proof}[Sketch]
    Define $E^n:=\underline U^n-\overline U^n$ and freeze a maximizing supersolution policy at each dynamic row, or an $\varepsilon$-maximizer followed by $\varepsilon\downarrow0$. Subtracting the two residual inequalities, we obtain
    \begin{equation*}
        (A^{\bar\alpha}E^{n-1})_i
        \leq\frac{E_i^n}{\Delta t}\leq0
    \end{equation*}
    on dynamic rows, with the same sign on data and copy rows. Since $(A^{\bar\alpha})^{-1}\geq0$ by \autoref{thm:constructed_operator_monotone}(ii), $E^n\leq0$ implies $E^{n-1}\leq0$. The terminal inequality and backward induction finish the proof. See \autoref{app:proof_upwind_comparison} for the details.
\end{proof}

\paragraph{Barles--Souganidis ingredients.}

The Barles--Souganidis route \cite{BarlesSouganidis1991} requires monotonicity, stability, and consistency of the same residual. Exact quadratic Laplacian moments and affine fan exactness give the smooth-test consistency of \autoref{prop:consistency}, proved in \autoref{app:gfdm_consistency}.

Rows exist when the pool sees enough of the cloud. For a pool radius above a fixed multiple of the fill distance, a node whose pool ball meets no boundary receives a row of type (P) by \autoref{cor:constructibility}, and a node whose ball meets the boundary in a chord receives a row of type (P) or (R) at a reflecting face and (P) or (D) at an exit window by \autoref{lem:mirror_constructibility}.

For a discrete solution $U=(U_i^n)$, write $u_{h,\Delta t}(t_n,x_i):=U_i^n$, with half-relaxed limits as $(h,\Delta t)\to0$ taken over these space--time grid points, and abbreviate $u_h:=u_{h,\Delta t}$ when the time step is understood. Along a refinement family, $U$ and the density history $M$ depend on $h$, which we suppress, and we write
\begin{equation*}
    C_H
    :=\sup_h\sup_{n,i}
    \abs*{H_{i,h}
    [0;M_i^n]}.
\end{equation*}

\paragraph{Standing hypotheses.}
Suppose $\partial\Omega$ is a finite union of $C^1$ arcs, and let $\Gamma_{\mathrm{c}}\subset\partial\Omega$ be the finite set of points at which $\partial\Omega$ is not $C^1$, together with the points of $\Gamma_{\mathrm{refl}}\cap\overline{\Gamma_{\mathrm{abs}}}$. The results of this section draw on the following.
\begin{enumerate}
    \item[(H1)] \emph{Data.} $g_T:=G(\cdot,m(T))\in C(\overline\Omega)$. The identity-row terminal data are consistent, in the sense that $\sup_{x_i\in\mcalX_h}\abs{g_{h,i}^{N_t}-g_T(x_i)}\to0$. The imposed absorbing-Dirichlet data are uniformly bounded, $g_{\mathrm{exit}}$ is bounded and continuous on $[0,T)\times\Gamma_{\mathrm{abs}}$, and $g_{h,i}^n\to g_{\mathrm{exit}}(t_0,x_0)$ along absorbing-Dirichlet rows with $(t_n,x_i)\to(t_0,x_0)\in[0,T)\times\Gamma_{\mathrm{abs}}$.
    \item[(H2)] \emph{Hamiltonian.} $b,\mcalL$ in \eqref{eq:bellman_hamiltonian} are continuous, and on bounded test gradients the active controls of the continuous and the discrete Hamiltonians lie in an $h$-independent compact set.
    \item[(H3)] \emph{Rows.} At each $h$, every interior node is assigned a row of type (P), (D) or (R) in the sense of \autoref{def:row_classification}, except that, for fixed $C_{\mathrm{v}},C_{\mathrm{cl}}>0$, a node within $C_{\mathrm{v}}h$ of $\Gamma_{\mathrm{c}}$ may instead take, in place of a refusal, the copy row \eqref{eq:copy_row} from a dynamic node within $C_{\mathrm{cl}}h$, carrying the drift of that node. These indices join $\mscrI_h^{\mathrm{cl}}$.
    \item[(H4)] \emph{Comparison.} For a prescribed density $m$, $H_m(t,x,p):=H(x,p,m(t,x))$ is convex in $p$, and there are constants $C_{\mathrm{q}}>0$ and $L\geq0$ and a modulus $\omega_H$ such that
        \begin{align}
            H_m(t,x,p)&\leq C_{\mathrm{q}}(1+\abs{p}^2),\label{eq:comparison_growth}\\
            H_m(t,x,q)-H_m(s,y,q)&\leq\omega_H(\abs{t-s}+\abs{x-y})+L(\abs{t-s}+\abs{x-y})\abs{q}.\label{eq:comparison_structure}
        \end{align}
\end{enumerate}
At a prescribed continuous density, \textup{(H2)} makes $C_H$ finite and bounds the drifts along $M$.

\begin{proposition}[Fixed-density structure and stability]\label{thm:viscosity_convergence}
    Let $d=2$. Assume \textup{(H3)} and $C_H<\infty$, and let $C_0$ bound the imposed terminal and absorbing-Dirichlet data uniformly in $h$. Then, at each fixed density iterate, the fully implicit
    residual $F_{h,\Delta t}$ is monotone in the Barles--Souganidis sense, and its solutions satisfy the constant-barrier bound
    $\abs{U_i^n}\leq C_0+C_HT$, uniform in $h$ and $\Delta t$, so they are uniformly $L^\infty$-stable.
\end{proposition}

\begin{proof}[Sketch]
    Under \textup{(H3)}, each fixed-policy residual is monotone, and the Bellman supremum preserves this property. Because the spatial rows annihilate constants, $C_0+C_H(T-t_n)$ and $-C_0-C_H(T-t_n)$ are discrete super- and subsolutions. \autoref{thm:upwind_comparison} gives the stated bound. See \autoref{app:proof_viscosity_convergence} for further details.
\end{proof}

Consistency, the remaining Barles--Souganidis ingredient, requires the density entering the Hamiltonian to converge, which holds when the density is prescribed.

\begin{theorem}[Convergence on polygons at a prescribed density]\label{thm:polygon}
    Let $\partial\Omega$ be a finite union of straight segments, and let each point of $\Gamma_{\mathrm{refl}}\cap\overline{\Gamma_{\mathrm{abs}}}$ lie in the relative interior of one of them. Let $m\in C([0,T]\times\overline\Omega)$ be prescribed, set $M_i^n:=m(t_n,x_i)$, and assume \textup{(H1)}, \textup{(H2)} and \textup{(H4)}. Let the family have local support, under the constants and the construction rules of \autoref{lem:mirror_constructibility}, with $\delta_i\geq(1+1/\sin(\pi/48))h_i$ at every node whose pool ball meets $\Gamma_{\mathrm{refl}}$. Then there are $C_{\mathrm{v}},C_{\mathrm{cl}}>0$ such that, with the copy rows that \textup{(H3)} allows near $\Gamma_{\mathrm{c}}$, for all sufficiently small $h$ every interior node at distance at least $C_{\mathrm{v}}h$ from $\Gamma_{\mathrm{c}}$ receives a row of type (P), (D) or (R), and discrete solutions satisfy
    \begin{equation*}
        u_{h,\Delta t}\to u
        \quad\text{locally uniformly on }[0,T)\times(\overline\Omega\setminus\Gamma_{\mathrm{c}})
        \quad\text{as }(h,\Delta t)\to0,
    \end{equation*}
    where $u$ is the unique bounded viscosity solution of the terminal--boundary-value problem on $Q_{\mathrm{c}}:=([0,T]\times(\overline\Omega\setminus\Gamma_{\mathrm{c}}))\setminus(\{T\}\times\overline{\Gamma_{\mathrm{abs}}})$. In particular, this solution exists.
\end{theorem}

\begin{proof}[Sketch]
    Away from $\Gamma_{\mathrm{c}}$, a ball of the pool radius meets the boundary in at most one chord, so \autoref{lem:mirror_constructibility} gives a row at every node there. Copy sources near $\Gamma_{\mathrm{c}}$ exist, which gives \textup{(H3)}. Let $\overline u$ and $\underline u$ be the upper and lower half-relaxed limits. In $\Omega$, \autoref{prop:consistency} and \autoref{thm:viscosity_convergence} make them viscosity subsolutions and supersolutions. At an upper contact of $\overline u-\varphi$ on a reflecting segment, a subsequence of type (P) rows gives $F[\varphi;m]\leq0$, and on a subsequence of type (R) rows \autoref{prop:mirror_boundary_expansion} gives
    \begin{equation*}
        F_{h,\Delta t}(\pi_h\varphi;M)
        =F[\varphi;m](t_n,x_i)
        +\frac{\sigma^2}{2}\frac{\widehat\gamma_{i,h}}{\delta_i}
        \partial_\nu\varphi(t_n,x_0)+O(1).
    \end{equation*}
    If $\partial_\nu\varphi(t_0,x_0)>0$, then $\sigma>0$ and \eqref{eq:bc_gamma_main} force this residual to $+\infty$, contradicting the subsolution inequality. So $\min\{F[\varphi;m],\partial_\nu\varphi\}\leq0$, which is the relaxed Neumann condition, and the reversed argument gives $\max\{F[\varphi;m],\partial_\nu\varphi\}\geq0$ at a lower contact of $\underline u-\varphi$. Exit windows are similar, with data rows in place of mirror rows. \autoref{lem:terminal_closure} gives the terminal inequalities, \autoref{prop:comparison} gives comparison off the vertices, and so $\overline u=\underline u$ on $Q_{\mathrm{c}}$. See \autoref{app:half_relaxed} and \autoref{app:polygon} for the details.
\end{proof}

\autoref{app:smooth_rate} proves a residual-to-value transfer estimate and, with it, a first-order rate on smooth data at a prescribed density.

\section{Numerical results}\label{sec:numerical_results}

\autoref{tab:benchmarks} lists the benchmarks and the part of the method each one tests. In the experiments, $N_{\mathrm{int}}$ is the number of interior nodes and $h$ denotes the nominal spacing $\sqrt{\abs{\Omega}/N_{\mathrm{int}}}$, not the fill distance of \autoref{def:stencil_geometry}. On reflecting domains the pool radius is about five nominal spacings, below the radius that \autoref{thm:polygon} asks at a reflecting face, so the signs, fan completeness and obliqueness are checked on each computed row. \autoref{app:numerical_provenance} gives the software, cloud generation, KDE bandwidths, row radii and seeding common to all runs, the stencil sizes and time steps of each benchmark, a cost comparison with a five-point scheme, and the census of row types and fans.

\begin{table}[htbp]
    \centering
    \footnotesize
    \caption{Benchmarks and the stage each one tests.}
    \label{tab:benchmarks}
    \begin{tabular}{lll}
        \toprule
        Benchmark & Boundary & Stage tested \\
        \midrule
        Floor plan & reflecting polygon & one HJB solve \\
        Linear--quadratic (LQ) & Dirichlet frame & one HJB solve \\
        Source-free & reflecting square & coupled, particles \\
        Evacuation & polygon, disc, exits & coupled, particles \\
        \bottomrule
    \end{tabular}
\end{table}

\subsection{Exact solutions on a floor plan}\label{sec:floor_plan}

We test \autoref{thm:polygon} on two rooms joined by a door, $\Omega=(0,5)\times(0,3)\setminus([2,3]\times[0,1]\cup[2,3]\times[2,3])$, with all walls reflecting and four re-entrant jambs. Take $H=\frac12\abs{p}^2-V(t,x)$ and $\sigma=\frac12$. The solutions
\begin{equation}\label{eq:floor_plan_solutions}
    u_1=-\sigma^2\log\parens[\big]{\tfrac32+e^{-\sigma^2\pi^2(T-t)}\cos\pi x\cos\pi y},
    \qquad
    u_2=f_5(x)f_3(y),
\end{equation}
with $f_L(s):=\frac1\pi\parens[\big]{\frac{\sin\pi s}{L\pi}-(1+\frac sL)\cos\pi s}$, so that $f_L'(s)=(1+s/L)\sin\pi s$, have $\partial_\nu u=0$ on every wall and a nonzero tangential derivative along it. The first solves the equation with $V\equiv0$ and $T=1$, and the second with $V$ its stationary residual and $T=5$. Their oscillations are $0.40$ and $0.77$. The first is even across every wall, which makes the fold of a mirror row exact, and the second is not.

Each scrambled Sobol cloud runs with $\Delta t=h$ and a pool of radius about $5.3h$ (\autoref{fig:exact_benchmarks}a).

The error on $u_1$ falls at a local order that peaks near $1.8$ and is $1.3$ at $N_{\mathrm{int}}=83200$. The error on $u_2$ stops falling between $N_{\mathrm{int}}=5200$ and $20800$, and then falls to $0.058$ and $0.045$ at $N_{\mathrm{int}}=46800$ and $83200$. Propagated through the scheme with the policy frozen at the exact solution, the residual of the exact solution reproduces this error to within $10\%$ from $N_{\mathrm{int}}=5200$, and it separates two parts of opposite sign. The fan part is first order, about $3.8h$. The Laplacian part, mostly on mirror rows, falls at about order $2.5$. They cancel near $N_{\mathrm{int}}=5200$, which makes the flat range, and beyond it the error follows the fan part. On $u_1$ the fan part is small, about $0.05h$, and the local order falls as the Laplacian part gives way to the first-order time step, about $0.19h$. Cutting $\Delta t$ fourfold on the cloud of $N_{\mathrm{int}}=20800$ halves the error. Between one and eight nodes per cloud take copy rows, all within $0.3$ pool radii of a vertex and none at a jamb.

\subsection{Rates and cost on a square}\label{sec:mms_coupled}

The decoupled LQ benchmark solves
\begin{equation*}
    -\partial_tu+\tfrac12\abs{\nabla u}^2-\tfrac{\sigma^2}{2}\Delta u=0,
    \qquad
    u(T,x)=G_s\abs{x-x_{\mathrm{exit}}}^2,
\end{equation*}
on $\Omega=(0,20)^2$ with $T=4$, $\sigma=0.3$, $G_s=\frac12$ and $x_{\mathrm{exit}}=(20,20)$. Its viscosity solution is the classical Riccati solution
\begin{equation*}
    u(t,x)=P(t)\abs{x-x_{\mathrm{exit}}}^2+\sigma^2\log\parens[\big]{1+2G_s(T-t)},
    \qquad
    P(t)=\parens[\big]{1/G_s+2(T-t)}^{-1},
\end{equation*}
which also gives the Dirichlet data on the frame. Because the Laplacian moments are exact on quadratics, this datum isolates the fan and time errors of the HJB stage.

The fan-upwind scheme converges at near first order on scattered Sobol and jittered clouds (\autoref{fig:exact_benchmarks}b). At matched degrees of freedom (DOF), the ratio of the meshfree error to that of a five-point finite-difference method (FDM) is $1.11$ to $1.20$ in $L^2$ and $0.57$ to $0.63$ in $L^\infty$, and the stencil construction makes the meshfree solve about twice as expensive in total (\autoref{app:complexity_table}).

\subsection{Full coupling on a reflecting square}\label{sec:full_coupling}

A manufactured solution would need a source that an empirical measure cannot transport. We therefore construct a smooth source-free exact MFG on the purely reflecting square $\Omega=(0,1)^2$, which tests the fully coupled computation, particle stage included. Let $T=1$, $D:=\sigma^2/2=1/8$, $\zeta=1/2$, $C=1$, and set
\begin{equation}\label{eq:source_free_reference}
    \begin{aligned}
        \omega&=4\pi, & A&=\frac{D\ln 5}{4},\\
        \phi(x)&=A\parens[\big]{\cos(\omega x_1)+\cos(\omega x_2)},\\
        m^*(x)&=Z^{-1}e^{-\phi(x)/D},
        & u^*(t,x)&=\phi(x)+C(T-t)^2,
    \end{aligned}
\end{equation}
where $Z=\int_\Omega e^{-\phi/D}\diff x$. Thus $m^*$ has unit mass and contrast $\max m^*/\min m^*=5$. For the quadratic-control running cost and Hamiltonian
\begin{equation}\label{eq:source_free_hamiltonian}
    \begin{aligned}
        \mcalL&=\tfrac12\abs{\alpha}^2+V(t,x)+\zeta m,\\
        H&=\tfrac12\abs{p}^2-V(t,x)-\zeta m,\\
        V(t,x)&=2C(T-t)-D\Delta\phi(x)+\tfrac12\abs{\nabla\phi(x)}^2-\zeta m^*(x),
    \end{aligned}
\end{equation}
take $m_0=m^*$ and $G=\phi$. Then $b^*=-\nabla\phi$, the HJB residual vanishes by the definition of $V$, and
\begin{equation*}
    D\nabla m^*-m^*b^*=0,\qquad
    \partial_\nu\phi=0\quad\hbox{on }\partial\Omega,
\end{equation*}
so $(u^*,m^*)$ solves the source-free reflected MFG exactly.

The exact law is stationary, and the coupling stays active away from the fixed point because the effective running potential contains $\zeta(m-m^*)$.

\begin{figure}[htbp]
    \centering
    \includegraphics[width=\textwidth]{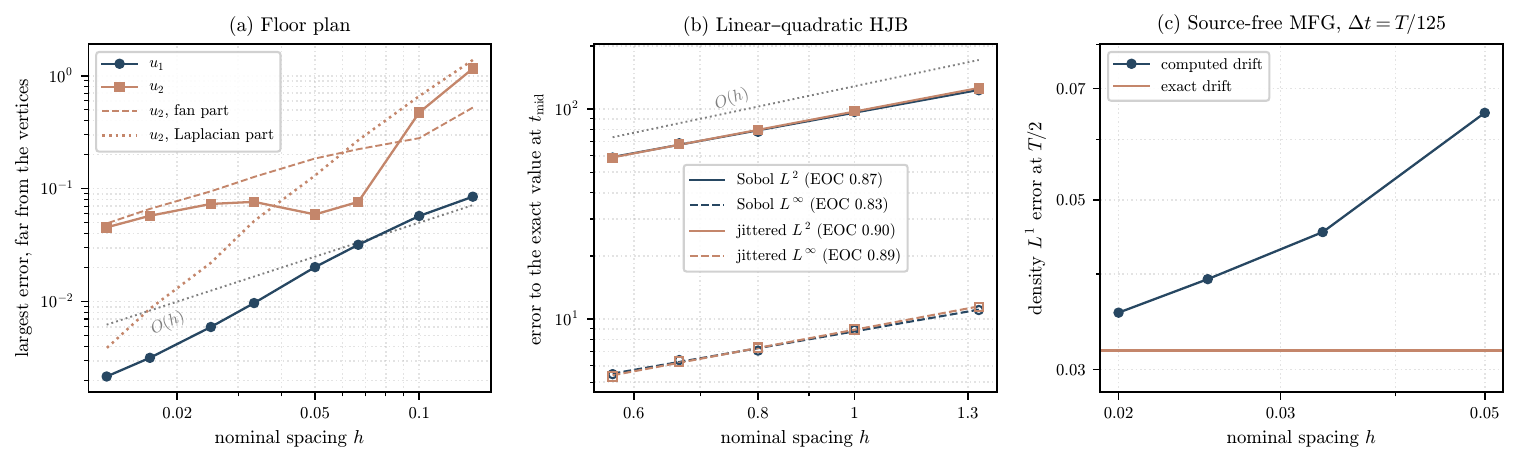}
    \caption{Errors against the exact solutions. (a) Floor plan, seed $42$: largest error over the time levels and over the nodes farther than one pool radius from every vertex, with the fan and Laplacian parts of the error on $u_2$. Up to $N_{\mathrm{int}}=20800$, seeds $1$ and $2$ differ from it by at most $21\%$. (b) Decoupled LQ HJB with $\Delta t=h/40$, on five Sobol and jittered refinements, with fitted slopes. The orders combine the $O(h)$ time-step and fan defects. (c) Source-free MFG, mean over seeds $42$, $1$ and $2$, at fixed $\Delta t$, $N_p$ and $h_{\mathrm{KDE}}$. The line is particles driven by the exact drift.}
    \label{fig:exact_benchmarks}
\end{figure}

At fixed $\Delta t=T/125$, $N_p=10^6$, and $h_{\mathrm{KDE}}=0.016$, refining the cloud from $N_{\mathrm{int}}=400$ to $2500$ reduces the three-seed mean mid-time density $L^1$ error from $0.0650$ to $0.0356$ (\autoref{fig:exact_benchmarks}c). Particles driven by the exact drift do not read the cloud and give $0.0318$.

\subsection{A non-convex game}\label{sec:main_result}

The evacuation game is posed on the domain $\Omega$ obtained from $[0,24]\times[0,12]$ by removing a circular obstacle of radius $r=1.4$ centered at $x_{\text{circle}}=(8,6)$ and a vertical door wall at $x=14$ of thickness $w=1.5$, punctured by a door gap of width $d_{\mathrm{door}}=1.5$ at $y\in[5.25,6.75]$. Two absorbing exits span the window $x\in[18,22]$ on the top and bottom outer walls, and the remaining walls reflect. The initial density is a Gaussian centered at $(1.5,6)$ with width $0.5$, so all mass starts in the left room and must route around the obstacle, pass the door, and commit to one exit. We take $\sigma=0.5$ and $T=40$. In this soft-wall game, wall and obstacle proximity enter the running cost through bounded exponential potentials $V_{\mathrm{wall}}$ and $V_{\mathrm{obs}}$, and congestion enters through $\zeta m$ with $\zeta=0.5$, so that
\begin{equation}\label{eq:main_result_H}
    H(x,p,m)=\tfrac{1}{2}\abs{p}^2-V_{\mathrm{wall}}(x)-V_{\mathrm{obs}}(x)-\zeta m,
\end{equation}
with the potentials
\begin{align*}
    V_{\mathrm{wall}}(x) &= C_{\mathrm{wall}}\parens[\big]{e^{-d_{\mathrm{obs}}^+(x) / \varepsilon_w} + e^{-d_{\text{outer}}(x) / \varepsilon_w}},
    & C_{\mathrm{wall}} &= 10,\ \varepsilon_w = 0.2,\\
    V_{\mathrm{obs}}(x) &= \kappa_{\mathrm{obs}} e^{-d_{\text{circle}}^+(x) / \lambda_{\mathrm{obs}}},
    & \kappa_{\mathrm{obs}} &= 30,\ \lambda_{\mathrm{obs}} = 3.0,
\end{align*}
where $d_{\mathrm{obs}}^+ = \max(d_{\mathrm{obs}}, 0)$ is the signed distance to the union of inner obstacles clamped nonnegative outside, $d_{\text{outer}}$ is the distance to the nearest outer-wall non-exit segment (the exit windows on the top and bottom walls are exempt from $V_{\mathrm{wall}}$), and $d_{\text{circle}}^+ = \max(\abs{x - x_{\text{circle}}} - r, 0)$ is the outward distance to the obstacle circle's boundary. The terminal cost is the geodesic soft-min squared-distance to the nearest exit:
\begin{equation}\label{eq:main_result_gT}
    g_T(x) = G_{\text{scale}} \cdot \operatorname{softmin}_\tau\parens[\big]{d_{\text{geo}}(x, \mathrm{exit}_{\mathrm{top}})^2, d_{\text{geo}}(x, \mathrm{exit}_{\mathrm{bot}})^2},
    \quad G_{\text{scale}} = 0.05,\ \tau = 4,
\end{equation}
computed on a $k$-nearest-neighbor (k-NN) graph with wall-crossing edges culled (so $d_{\text{geo}}$ respects obstacle topology). Here $\operatorname{softmin}_\tau(s_1,s_2) := -\tau \log(e^{-s_1/\tau} + e^{-s_2/\tau})$ smooths the medial-axis kink between the two exit basins at $y = 6$, and the exit datum is $g_{\mathrm{exit}}\equiv0$.

\begin{figure}[htbp]
    \centering
    \includegraphics[width=0.9\textwidth]{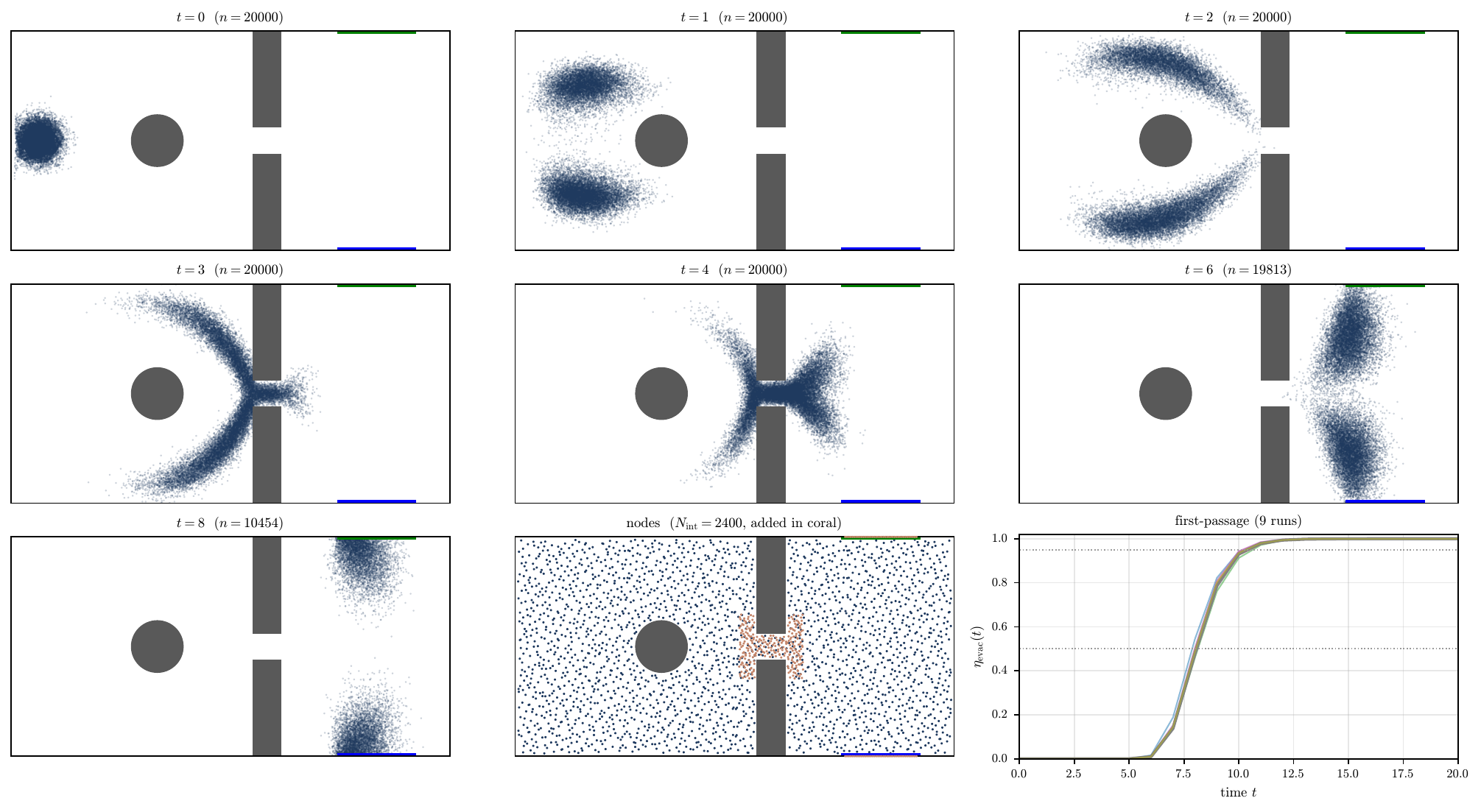}
    \caption{Evacuation (seed 42, $N_{\mathrm{int}}=2400$). Particles split around the obstacle and reach both exits. The eighth panel shows the nodes, with the added ones in coral. First-passage curves cover all nine runs.}
    \label{fig:density_evolution}
\end{figure}

We apply our construction to this game (\autoref{fig:density_evolution}). Each cloud carries $280$ further nodes around the door, and the exit nodes are spaced $0.1$ apart. Across nine clouds and $16{,}920$ interior rows, $1684$ are mirror rows and $41$ next to the boundary vertices are copy rows, and every implicit operator satisfies the M-matrix sign conditions to solver tolerance (\autoref{app:row_census}). The median of the first-passage time is $8$ to $9$ and its $95$th percentile is $11$ on all nine clouds, where the snapshots are one time unit apart. The top-exit fraction lies between $0.41$ and $0.58$ on the nine clouds, against $0.21$ to $0.93$ without the added nodes.

\section{Conclusion}\label{sec:conclusion}

We developed a meshfree--particle method for mean field games. It solves the HJB equation on a fixed scattered cloud, carries the Fokker--Planck law by reflected and killed particles, and closes the fictitious-play loop by interpolation and mass-preserving kernel density estimation. Monotonicity is imposed and checked row by row, so the discrete comparison principle holds on a scattered cloud at any cell-P\'eclet number, while on a mesh a single violating edge voids the angle-condition guarantee for the whole matrix. On polygons, re-entrant corners included, one HJB solve at a prescribed density converges to the viscosity solution off the boundary vertices. Exact solutions on a floor plan and on a square test the rows of the theorem and measure their rates, the source-free game measures the error of the fully coupled computation, and on the evacuation game every implicit operator on nine clouds satisfies the M-matrix sign conditions to solver tolerance.

What remains open is the coupled limit. \autoref{thm:polygon} holds at a prescribed density. The convergence of the damped random recursion of \autoref{alg:Hybrid-Eulerian-Lagrangian}, and of its density iterates, to a mean field equilibrium is not analyzed. A particle forward stage preserves mass exactly but is not the discrete adjoint of the backward operator, so the route to coupled rates through one stability constant \cite{Osborne2026Rates} does not apply. The coupled runs use a fixed number of outer iterations, because a tolerance on the fixed-point residual is not meaningful below the sampling noise of one best response. The source-free instance, the only exact coupled benchmark, is stationary, so the particle stage is tested under full coupling only against a stationary law. A mean-field-coupled LQ instance, whose law is Gaussian with an explicit mean, is the natural next test.

Hypothesis \textup{(H3)} is uniform in $h$ and cannot follow from finitely many clouds. \autoref{thm:polygon} obtains it on polygons for one ray per cell, with $\delta_i\geq(1+1/\sin(\pi/8))h$ at all nodes and $\delta_i\geq(1+1/\sin(\pi/48))h_i$ where the pool ball meets a reflecting face. The computed pools fall below the second bound and take their fans from a search over the pool, not from the cells. On the evacuation domain, \autoref{prop:comparison} gives comparison at a prescribed density off the boundary vertices, and \autoref{thm:viscosity_limit}\textup{(ii)} covers planar reflecting segments and exit windows away from junctions. The circle lies outside it, since the mirror rows there reflect across a tangent line. The terminal data are graph distances from a cloud-dependent source set, and \textup{(H1)} is not established for them. Rates for nonsmooth solutions on scattered clouds would require Krylov-type estimates.

\printbibliography[heading=bibintoc,title={References}]

\appendix

\section{Discrete constructions}\label{app:discrete_details}

\subsection{Constructibility}\label{app:constructibility}

\begin{definition}[Visible availability and coverage]\label{def:visibly_available}
    A set $S\subset\overline\Omega$ is \emph{visibly available from $x_i$} if each edge $[x_i,y]$ with $y\in S$ is visible. Given a family of angular sectors about $x_i$, the node has \emph{visible coverage} of that family if each sector contains a closed ball of radius $h_i$ that is visibly available from $x_i$.
\end{definition}

\begin{lemma}[Fan from annular cells]\label{lem:fan_cells}
    Let $d=2$, let $N_{\mathrm{sec}}\geq4$ be an integer, let $\theta_0:=2\pi/N_{\mathrm{sec}}$, let $0<\omega<\theta_0$, and let $0<c_r<C_r$ satisfy
    \begin{equation}\label{eq:fan_cell_geometry}
        \frac{C_r-c_r}{2}\geq1,
        \qquad
        \frac{C_r+c_r}{2}\sin\frac{\omega}{2}\geq1.
    \end{equation}
    The smallest such $C_r$ is $1+1/\sin(\omega/2)$, with $c_r=C_r-2$. For $j=0,\ldots,N_{\mathrm{sec}}-1$, define the annular cell
    \begin{align*}
        S_{i,j}
        &:=\bigl\{x:c_rh_i\leq\abs{x-x_i}\leq C_rh_i,\\
        \arg(x-x_i)&\in[(j+\tfrac12)\theta_0-\omega/2,(j+\tfrac12)\theta_0+\omega/2]\bigr\}.
    \end{align*}
    If $x_i$ has visible coverage of $\{S_{i,j}\}_j$ (\autoref{def:visibly_available}), $\delta_i\geq C_rh_i$ and the search radius is at least $C_rh_i$, one ray selected from each $S_{i,j}$ forms a complete fan for any cap $\theta_+\geq\theta_0+\omega$, with
    \begin{equation*}
        \theta_0-\omega
        \leq\angle(z_{i,j},z_{i,j+1})
        \leq\theta_0+\omega<\pi.
    \end{equation*}
\end{lemma}

\begin{proof}
    The point at radius $r_0=(C_r+c_r)h_i/2$ on the angular bisector of $S_{i,j}$ has distance $(C_r-c_r)h_i/2$ from the radial faces and distance $r_0\sin(\omega/2)$ from the angular faces. Condition \eqref{eq:fan_cell_geometry} therefore places a closed ball of radius $h_i$ inside each cell, so the coverage hypothesis is not vacuous. Visible coverage makes each edge from $x_i$ to a point in such a ball visible in the sense of \autoref{def:visibility}, and the search-radius hypothesis puts the annulus inside the admissible search region. The interior-ball condition $B(x_i,C_rh_i)\subset\Omega$ gives visible coverage and is the unoccluded special case.

    Apply the fill-distance definition at the center of each inscribed ball, which lies in $\Omega$ because $\Omega$ is the interior of $\overline\Omega$. A cloud point lies within distance $h_i$ and hence remains inside that annular cell. Selecting one such point per cell gives rays with $c_r\leq\abs{z_{i,j}}/h_i\leq C_r$ and adjacent angles in $[\theta_0-\omega,\theta_0+\omega]\Subset(0,\pi)$, and their cyclic order makes the adjacent cones tile $\mbbR^2$.
\end{proof}

\begin{corollary}[Constructibility from visible coverage]\label{cor:constructibility}
    Let $d=2$, let the refinement family have local support, and let $0<c_r<C_r$ satisfy \eqref{eq:fan_cell_geometry} at $\omega=\pi/4$ and at the $\omega$ of \autoref{lem:fan_cells}.
    \begin{enumerate}
        \item[(i)] Let $\delta_i\geq C_rh_i$, and let every sector
            $\{c_rh_i\leq\abs{x-x_i}\leq C_rh_i,\ \arg(x-x_i)\in[\varphi,\varphi+\pi/4]\}$,
            $\varphi\in[0,2\pi)$, contain a closed ball of radius $h_i$ visibly available
            from $x_i$. Then \eqref{eq:qp_problem} is feasible at $x_i$.
        \item[(ii)] If moreover $c_\delta\geq C_r$ and $\theta_+\geq\theta_0+\omega$, then, for every $K\Subset\Omega$ and all sufficiently small $h$, every node of $K$ receives a row of type (P) whose fan takes one ray per annular cell of \autoref{lem:fan_cells}.
    \end{enumerate}
\end{corollary}

\begin{proof}
    The two inequalities place a radius-$h_i$ ball inside each $\pi/4$ annular sector. The first controls the radial faces, and the second controls the angular faces at the mid-radius. Visible coverage makes the points of each such ball visible from $x_i$, so the fill distance at its center places a cloud point in it and every resulting edge is visible (\autoref{def:visibility}). Since the sectors are taken at all orientations, no $\pi/4$ cone about $x_i$ is empty, which is the hypothesis of Seibold's cone criterion (the cone of \cite[Thm.~9]{Seibold2008Minimal} has a direction depending on the moment data, and the sufficiency is quantified over that direction). It then provides a nonnegative stencil satisfying the constant, linear, and quadratic Laplacian moments.

    For the uniform statement, $h_i\leq h$ by the definition of $h$, so \autoref{def:local_support} and $c_\delta\geq C_r$ give $\delta_i\geq c_\delta h\geq C_rh_i$. With $h_K:=\operatorname{dist}(K,\partial\Omega)/C_r$, every $x_i\in\mcalX_h\cap K$ has $B(x_i,C_rh_i)\subset\Omega$ once $h<h_K$. The coverage and the fan cells are then available at all nodes of $K$, so each receives a row of type (P), with the fan of \autoref{lem:fan_cells}.
\end{proof}

If $B(x_i,C_rh_i)\subset\Omega$, nothing is occluded and the coverage hypothesis of \textup{(i)} holds. The lower bound on the pool radius in \autoref{def:local_support} serves the hypothesis $\delta_i\geq C_rh_i$. Without it, a family meeting only $\delta_i\leq C_\delta h$ may have $\delta_i<C_rh_i$ at all nodes and no complete fan anywhere. The corollary and \autoref{lem:fan_cells} together give the PDE branch at any node with visible coverage, $\delta_i\geq C_rh_i$ and search radius at least $C_rh_i$, occluded or not, since (P) takes precedence in \autoref{def:row_classification}.

\subsection{Mirror constructibility on a planar face}\label{app:mirror_constructibility}

\begin{lemma}[Mirror constructibility on a planar face]\label{lem:mirror_constructibility}
    Let $d=2$, let the family have local support, let $C_r:=1+1/\sin(\pi/8)$ and $c_r:=C_r-2$, and let the $\omega$ of \autoref{lem:fan_cells} be at least $\pi/4$. Let $c_\delta\geq C_r$, $\kappa\geq C_r/c_\delta$ and $\theta_+\geq\theta_0+\omega$, and let a physical fan take one ray per annular cell of \autoref{lem:fan_cells}. Let $x_i$ be an interior node such that $B(x_i,\delta_i)\cap\partial\Omega$ is empty or is a chord of the ball lying in $\Gamma_{\mathrm{refl}}$ or in $\Gamma_{\mathrm{abs}}$.
    \begin{enumerate}
        \item[(a)] If the intersection is empty, then $x_i$ receives a row of type (P).
        \item[(b)] If the chord lies in $\Gamma_{\mathrm{abs}}$, then $x_i$ receives a row of type (P) or (D).
        \item[(c)] Let the chord lie in $\Gamma_{\mathrm{refl}}$, and let $\delta_i\geq C'h_i$, where $C':=1+1/\sin(\pi/48)$. Let the mirror program take its ghosts from $\mcalV_i\setminus\{x_i\}$ across the nearest face, with $\gamma_0\in(0,1]$ in \eqref{eq:bc_gamma_main}, and let a reflected fan take one ray per annular cell. Then $x_i$ receives a row of type (P) or (R).
    \end{enumerate}
\end{lemma}

\begin{proof}
    If the intersection is empty, then $B(x_i,C_rh_i)\subset\Omega$, and the step \emph{Nodes away from the face} below gives (a). Otherwise use tangential--normal coordinates $x=(y,r)$, with the chord on $\{r=0\}$, physical side $\{r>0\}$ and reflection $R(y,r):=(y,-r)$, and write $s_{ii}>0$ for the inward distance of the center. The chord holds the only boundary points of the ball, so its face is the nearest to $x_i$, and (D) does not occur when the chord lies in $\Gamma_{\mathrm{refl}}$.

    \emph{Coverage by the augmented pool.} The chord divides the ball into two connected sets that do not meet $\partial\Omega$. The one containing $x_i$ lies in $\Omega$, and the other lies outside $\overline\Omega$ because $\Omega$ lies locally on one side of its boundary. Thus $D_i:=B(x_i,\delta_i)\cap\Omega=B(x_i,\delta_i)\cap\{r>0\}$ is convex, so all its points are visible from $x_i$, and \autoref{def:stencil_geometry} places a pool node within $h_i$ of each point of $\overline{D_i}$. For $p=(y,r)$ with $r\leq0$ we have
    \begin{equation}\label{eq:reflection_closer}
        \abs{R(p)-x_i}^2=\abs{y-y_i}^2+(\abs{r}-s_{ii})^2\leq\abs{y-y_i}^2+(\abs{r}+s_{ii})^2=\abs{p-x_i}^2.
    \end{equation}
    Hence $R$ maps $B(x_i,\delta_i)\cap\{r\leq0\}$ into $\overline{D_i}$. Let $B$ be a closed ball of radius $h_i$ with center $c\in B(x_i,\delta_i)$. If $c$ lies in $\{r\geq0\}$, a pool node lies in $B$. Otherwise a pool node lies within $h_i$ of $R(c)$, and its ghost lies in $B$. So every such ball contains a pool node or the ghost of a pool node.

    \emph{Nodes away from the face.} If $s_{ii}\geq C_rh_i$, then $B(x_i,C_rh_i)\subset\Omega$. Since $\delta_i\geq c_\delta h\geq C_rh_i$, \autoref{cor:constructibility}(i) and \autoref{lem:fan_cells} give a row of type (P).

    \emph{Nodes near the face: a row on selected points.} Let $s_{ii}<C_rh_i$, so that $s_{ii}<(C_r/c_\delta)\delta_i\leq\kappa\delta_i$ and the face is within reach. Write $n$ for the $r$-component of a displacement from $x_i$, and let $S$ be the set of displacements to pool nodes with $n\geq0$, other than $x_i$, and to ghosts of pool nodes $x_j$ with $s_{ij}>s_{ii}$, all of length at most $C'h_i$. We show that $S$ meets every closed $\pi/4$ sector about the origin.

    Reversing the tangential axis if necessary, let the sector run from elevation $-\beta_2$ to elevation $\beta_1$ above the tangential direction, with $\beta_1+\beta_2=\pi/4$ and $\beta_1,\beta_2\leq5\pi/8$, and put $\rho:=1/\sin(\pi/48)$. A closed ball of radius $h_i$ with center at distance $\rho h_i$ from $x_i$ has its center in $B(x_i,\delta_i)$ and lies in $\{\abs{z}\leq C'h_i\}$. If $\beta_1\geq\pi/24$, center it on the bisector of the subsector from elevation $\beta_1-\pi/24$ to $\beta_1$. It lies in that subsector, hence in $\{n\geq0\}$, and the coverage step gives a pool node in it, which is not $x_i$ because the ball misses the origin. If $\beta_1<\pi/24$, then $\beta_2>5\pi/24$. Center the ball at elevation $-\psi$ with $\psi:=\beta_2-\pi/48$, so that it is tangent to the lower face and lies in the sector. Since $9\pi/48<\psi\leq29\pi/48$, its points satisfy
    \begin{equation*}
        n\leq-\rho h_i\sin\psi+h_i
        \leq-\parens[\Big]{\frac{\sin(3\pi/16)}{\sin(\pi/48)}-1}h_i
        <-2C_rh_i
        <-2s_{ii}.
    \end{equation*}
    The center has $r<-s_{ii}<0$, so the coverage step gives the ghost of a pool node $x_j$ in the ball, and that ghost has $r<-s_{ii}$, that is, $s_{ij}>s_{ii}$.

    We claim that there are $w_k\geq0$ with $\sum_kw_kz_k=0$ and $\sum_kw_kz_kz_k^{T}=2\Id$, the sums running over $z_k\in S$. Otherwise, Farkas' lemma gives $a\in\mbbR^2$ and a symmetric matrix $A$ with $\operatorname{tr}A<0$ and $f(z):=a\cdot z+z^{T}Az\geq0$ on $S$. Rotate so that $A=\operatorname{diag}(\lambda_1,\lambda_2)$ with $\lambda_2\leq\lambda_1$, whence $\lambda_2<0$, and choose signs $\epsilon_\ell$ with $\epsilon_\ell a_\ell\leq0$. On the closed $\pi/4$ sector $\{\epsilon_1z_1\geq0,\ \epsilon_2z_2\geq0,\ \abs{z_1}\leq\abs{z_2}\}$ minus the origin we have $z_2\neq0$ and $f(z)\leq\lambda_1z_1^2+\lambda_2z_2^2\leq\max(\lambda_2,\lambda_1+\lambda_2)z_2^2<0$. Since $S$ meets this sector, this is a contradiction. The weights $\widetilde L_{ik}:=w_k$, with the center weight fixed by the constant moment, form a nonnegative row supported on $S$ that is exact on $\msfP_2$.

    \emph{Ghost moment.} Write $n_k$ for the $r$-component of $x_k-x_i$ for this row. Exactness on $\msfP_2$ gives
    \begin{equation*}
        \sum_{k\neq i}\widetilde L_{ik}n_k=0,
        \qquad
        \sum_{k\neq i}\widetilde L_{ik}n_k^2=2.
    \end{equation*}
    These identities and $\abs{n_k}\leq C'h_i$ give
    \begin{equation*}
        \sum_{n_k<0}\widetilde L_{ik}\abs{n_k}
        =\frac12\sum_{k\neq i}\widetilde L_{ik}\abs{n_k}
        \geq\frac{1}{2C'h_i}\sum_{k\neq i}\widetilde L_{ik}n_k^2
        =\frac{1}{C'h_i}.
    \end{equation*}
    A point of $S$ with $n_k<0$ is the ghost of some $x_j$ with $s_{ij}>s_{ii}$, so $\abs{n_k}=s_{ij}+s_{ii}<2s_{ij}$. Therefore
    \begin{equation*}
        \sum_{j\in\mscrJ_i^{\mathrm{ghost}}}\widetilde L^{\mathrm{ghost}}_{ij}s_{ij}
        \geq\frac12\sum_{n_k<0}\widetilde L_{ik}\abs{n_k}
        \geq\frac{1}{2C'h_i}.
    \end{equation*}
    Since $\delta_i\geq C'h_i$, this row has $\widehat\gamma_{i,h}\geq1\geq\gamma_0$ in \eqref{eq:bc_gamma_main}. The mirror program is therefore feasible, and its minimizer satisfies \eqref{eq:bc_gamma_main} by its constraint.

    \emph{Folded fan.} An annular cell contains a closed ball of radius $h_i$ by \eqref{eq:fan_cell_geometry}, which lies in $B(x_i,\delta_i)$. By the coverage step it contains a pool node or a ghost, this time possibly the ghost of the center, whose ray folds onto the diagonal and contributes nothing to the row. One point per cell gives a complete fan by \autoref{lem:fan_cells}, and folding keeps its coefficients nonnegative.

    \emph{Conclusion.} The node now has a row of type (P) or (R), and (P) takes precedence.

    \emph{Exit windows.} If the chord lies in $\Gamma_{\mathrm{abs}}$, a node at inward distance at least $C_rh_i$ receives a row of type (P) as above, and a nearer node has the exit within reach, so it receives a row of type (P) or (D).
\end{proof}

\subsection{Howard policy iteration}\label{app:howard_solver}

For each nodal vector $W$, the discrete Hamiltonian of the quadratic instance retains the running cost in Bellman form:
\begin{equation}\label{eq:quadratic_discrete_hamiltonian}
    \begin{aligned}
        H_{i,h}[W;m_i]
        &=\sup_{\alpha\in\mbbR^2}
        \ensemble*{-(A^\alpha_{\mathrm{fan}}W)_i
        -\mcalL(x_i,m_i,\alpha)}\\
        &=-\min_{\alpha\in\mbbR^2}
        \ensemble*{(A^\alpha_{\mathrm{fan}}W)_i
        +\mcalL(x_i,m_i,\alpha)}.
    \end{aligned}
\end{equation}
At each backward step, Howard iteration alternates exact policy improvement in \eqref{eq:quadratic_discrete_hamiltonian} with a policy-evaluation linear solve. These steps are
\begin{equation}\label{eq:howard_alpha}
    \alpha^{(s+1)}_i
    \in\argmin_{\alpha\in\mbbR^2}
    \ensemble*{(A^\alpha_{\mathrm{fan}}U^{(s),n-1})_i
    +\mcalL\parens*{x_i,M_i^{(k),n},\alpha}},
    \qquad i=1,\ldots,N_{\mathrm{col}},
\end{equation}
and
\begin{equation}\label{eq:howard_linsolve}
    \parens[\big]{A^{\alpha^{(s+1)}}U^{(s+1),n-1}}_i
    =\frac{U_i^{(k+1),n}}{\Delta t}
    +\mcalL\parens*{x_i,M_i^{(k),n},\alpha_i^{(s+1)}}
\end{equation}
on the dynamic rows, with $A^\alpha$ the matrix \eqref{eq:implicit_matrix}.

Exact-policy Howard iteration is well defined and monotone, and converges to the discrete fixed point, when the control set is compact and the policy-evaluation matrices are inverse-positive, continuous in the policy and row-separable \cite[Thm.~2.1]{BokanowskiMarosoZidani2009}. That theorem is stated in minimization form and applies here through $U\mapsto-U$. The iterates are discrete supersolutions, so they decrease and stay above the constant barrier of \autoref{thm:viscosity_convergence}. At fixed $h$ their discrete gradients are therefore bounded. When $\mcalA$ is not compact, \autoref{lem:fan_existence} gives $\abs{(A^\alpha_{\mathrm{fan}}W)_i}\leq C_\theta\abs{b(x_i,m_i,\alpha)}\max_j\abs{W_j-W_i}/\min_j\abs{z_{ij}}$, which grows at most linearly in $\abs{b}$, while $\mcalL/(1+\abs{b})\to\infty$ as $\abs{\alpha}\to\infty$. Along the iterates the maximizing controls therefore lie in one compact set, and the theorem applies with $\mcalA$ replaced by it. \autoref{thm:constructed_operator_monotone}\textup{(ii)} gives inverse positivity on rows of types (P) and (R). Row $i$ of $A^\alpha$ depends on $\alpha_i$ alone because it contains the fan coefficients $\lambda_{ij}^{\alpha_i}$ and the policy-independent Laplacian row. On a complete fan these coefficients depend continuously on the drift, which is continuous in the control by \textup{(H2)}, so $A^\alpha$ is continuous in the policy. Hence the discrete system \eqref{eq:discrete_hjb_system} has a solution at each $h$ and $\Delta t$, unique by \autoref{thm:upwind_comparison}.

For a retained sector $(j,j+1)$, write $C_{i,j}:=\operatorname{cone}(z_{i,j},z_{i,j+1})$, and collect the two rays in $Z_{i,j}:=[\,z_{i,j}\ \ z_{i,j+1}\,]$ and the two nodal differences in $\Delta U_{i,j}:=(U_j-U_i,\,U_{j+1}-U_i)^\top$. The \emph{sector gradient}
\begin{equation*}
    g_{i,j}:=Z_{i,j}^{-\top}\Delta U_{i,j}
\end{equation*}
is the unique vector with $(A^\alpha_{\mathrm{fan}}U)_i=\alpha\cdot g_{i,j}$ for all $\alpha\in C_{i,j}$, where $\alpha=Z_{i,j}\lambda$ and $\lambda\geq0$.

For the quadratic instance, the sector objective is linear in $\alpha$ through $(A^\alpha_{\mathrm{fan}}U)_i=\alpha\cdot g_{i,j}$, so we have
\begin{equation*}
    \alpha_{i,j}^{*}
    =\operatorname{proj}_{C_{i,j}}
    \parens*{-\frac{g_{i,j}}{c}},
\end{equation*}
and finite comparison over the retained sectors amounts to the exact fan policy in \eqref{eq:howard_alpha}.

\subsection{KDE operator}\label{app:kde_construction}

The domain kernel of \eqref{eq:kde_operator} is built by the method of images. With $K$ the standard Gaussian, $h_{\mathrm{KDE}}$ a single absolute length, and $\mscrG$ the group generated by the reflections in the planes of $\Gamma_{\mathrm{refl}}$,
\begin{equation}\label{eq:image_kernel}
    K^\Omega_{h_{\mathrm{KDE}}}(x,y)
    =\sum_{g\in\mscrG}K_{h_{\mathrm{KDE}}}\parens[\big]{x-g(y)},
    \qquad
    K_{h_{\mathrm{KDE}}}(x):=h_{\mathrm{KDE}}^{-d}K(x/h_{\mathrm{KDE}}).
\end{equation}
The normalization in \eqref{eq:kde_operator} holds exactly whenever the images of $\Omega$ under $\mscrG$ tile $\mbbR^d$. On the rectangle $\Omega=\prod_{k=1}^d(0,L_k)$ the sum is a product of one-dimensional Neumann heat kernels at time $h_{\mathrm{KDE}}^2/2$, each with the cosine series of \cite{Botev2010}. The \emph{spectral realization} evaluates that series by a discrete cosine transform of the binned empirical measure and truncates no image, so it computes \eqref{eq:kde_operator} itself. The \emph{ghost-image reconstruction} serves geometries with no closed form. It keeps only the reflections in faces within a fixed multiple of $h_{\mathrm{KDE}}$, takes no images across interior obstacles, and truncates the Gaussian at radius $4h_{\mathrm{KDE}}$. Its kernel $\widetilde K$ therefore loses a source-dependent fraction of each deposit, largest near a wall, and the reconstruction is rescaled to the input mass,
\begin{equation}\label{eq:ghost_rescaled}
    \mcalK^{\mathrm{gh}}_{h_{\mathrm{KDE}}}[\mu]
    :=\frac{\mu(\overline\Omega)}{\int_{\overline\Omega}\widetilde{\mcalK}[\mu]}\,\widetilde{\mcalK}[\mu],
    \qquad
    \widetilde{\mcalK}[\mu](x):=\int_{\overline\Omega}\widetilde K(x,y)\diff\mu(y),
\end{equation}
with $\mcalK^{\mathrm{gh}}_{h_{\mathrm{KDE}}}[0]:=0$. The quotient makes $\mcalK^{\mathrm{gh}}_{h_{\mathrm{KDE}}}$ positively homogeneous but not additive, hence nonlinear, and it spreads the restored mass over the whole domain. Both realizations preserve the input mass, so neither restores mass killed at an exit. Consistency of kernel density estimation is classical \cite{Silverman1986}, and the boundary bias of reflected kernels is analyzed in \cite{Botev2010}.

\section{Proofs for the HJB analysis}\label{app:hjb_proofs}

\subsection{Proof of \autoref{lem:fan_existence} (simplicial-fan representation)}\label{app:proof_fan_existence}

Fix $v\in\mbbR^2\setminus\{0\}$. Because the rays form a complete cyclic fan with all adjacent angles in $(0,\pi)$, $v$ belongs to one adjacent two-ray cone. For that pair the system $[\,z_{i,j}\ z_{i,j+1}\,]\lambda=v$ has $\lambda\geq0$. Write $a_k:=\lambda_k\abs{z_k}$ for the two rays, at angle $\theta\leq\theta_+$. Since $a_1a_2\leq(a_1+a_2)^2/4$,
\begin{equation*}
    \abs{v}^2=(a_1+a_2)^2-2a_1a_2(1-\cos\theta)\geq(a_1+a_2)^2\cos^2(\theta/2),
\end{equation*}
so $a_1+a_2\leq\abs{v}/\cos(\theta/2)\leq C_\theta\abs{v}$, with equality when $a_1=a_2$ and $\theta=\theta_+$.

\subsection{Proof of \autoref{thm:constructed_operator_monotone} (algebraic monotonicity)}\label{app:proof_constructed_operator}

    \emph{Base and enlarged positive rows.} Feasibility of \eqref{eq:qp_problem} gives $L_{ij} \geq 0$ $(j \neq i)$ from its constraint, and $q\equiv1$ in its exactness condition gives $\sum_j L_{ij}=0$. The enlarged construction solves the same program on a larger visibility-safe stencil (\autoref{def:visibility}), with identical sign structure and zero row sum.

    \emph{Mirror-ghost rows.} The augmented program returns $\widetilde L_{ik}\geq0$ $(k\neq i)$ and $\sum_k\widetilde L_{ik}=0$. Folding adds each ghost weight to the weight of its physical source. For an off-center source this sums nonnegative weights, and a ghost of the center folds onto the diagonal. Each augmented weight is collected exactly once, whatever the number of faces, so the folded row keeps the sign structure and the zero row sum. The same collection of terms folds each guarded ghost ray onto its physical source in the fan.

    \emph{Assembly and (ii).} Every assembled non-Dirichlet Laplacian row is thus Metzler with zero row sum. The simplicial fan adds off-diagonals $\lambda_{ij} \geq 0$ acting through the differences $u_j - u_i$, hence zero advective row sums, proving (i) for $\Delta_h + A_{\mathrm{fan}}^\alpha$. On these dynamic rows $A^\alpha$ of \eqref{eq:implicit_matrix} has off-diagonals $-\frac{\sigma^2}{2}L_{ij}-\lambda_{ij}\leq0$ and row sums $1/\Delta t$, so each such row is strictly diagonally dominant. Absorbing-Dirichlet rows are identity rows. Copy rows \eqref{eq:copy_row} are only weakly diagonally dominant, so we eliminate the copied unknowns. Ordering the non-copy rows first,
    \begin{equation*}
        A^\alpha=
        \begin{pmatrix}A_{11}&A_{12}\\-C&\Id
        \end{pmatrix},
        \qquad
        (A^\alpha)^{-1}=
        \begin{pmatrix}\widehat A^{-1}&-\widehat A^{-1}A_{12}\\C\widehat A^{-1}&\Id-C\widehat A^{-1}A_{12}
        \end{pmatrix}.
    \end{equation*}
    Here $A_{11}$ and $A_{12}\leq0$ hold the entries of the dynamic and data rows on non-copied and copied nodes, and $C\geq0$ has one unit entry per copy row. Thus $\widehat A=A_{11}+A_{12}C$ has nonpositive off-diagonals. Each row of $C$ sums to one, so each row of $\widehat A$ has the row sum of the corresponding row of $A^\alpha$, namely $1/\Delta t$ on dynamic rows and $1$ on data rows. Hence $\widehat A$ is strictly diagonally dominant, and so a nonsingular M-matrix with $\widehat A^{-1}\geq0$ \cite[Chap.~6, Thm.~2.3]{BermanPlemmons1994}. Block multiplication verifies the displayed inverse, which is nonnegative because $A_{12}\leq0$ and $C\geq0$. As $A^\alpha$ is a $Z$-matrix with a nonnegative inverse, it is a nonsingular M-matrix \cite[Chap.~6]{BermanPlemmons1994} for all $\Delta t>0$ and all fixed admissible policies.

\subsection{Proof of \autoref{thm:upwind_comparison} (P\'eclet-free discrete comparison principle)}\label{app:proof_upwind_comparison}

    For a fixed rowwise policy $\alpha$, let $A^\alpha$ be the matrix \eqref{eq:implicit_matrix}. The hypotheses are precisely the algebraic inputs of the block argument in \autoref{thm:constructed_operator_monotone}(ii), so $(A^\alpha)^{-1}\geq0$.

    Set $E^n:=\underline U^n-\overline U^n$. The terminal data and copy inequalities give $E^{N_t}\leq0$. Suppose $E^n\leq0$, fix $\varepsilon\geq0$, and at each dynamic row choose an $\varepsilon$-maximizer $\bar\alpha_i$ of the supremum for $\overline U^{n-1}$, with $\varepsilon=0$ when it is attained. Evaluating the subsolution Hamiltonian at the same policy and subtracting the supersolution inequality, we obtain
    \begin{equation*}
        (A^{\bar\alpha}E^{n-1})_i
        \leq\frac{E_i^n}{\Delta t}+\varepsilon\leq\varepsilon.
    \end{equation*}
    On data rows $E_i^{n-1}\leq0$, while on copied nodes
    $E_i^{n-1}-E_{\iota(i)}^{n-1}\leq0$. The constant vector $\Delta t\,\mathbf 1$ has $A^{\bar\alpha}(\Delta t\,\mathbf 1)$ equal to $1$ on dynamic rows, $\Delta t$ on data rows and $0$ on copy rows, so $A^{\bar\alpha}(E^{n-1}-\varepsilon\Delta t\,\mathbf 1)\leq0$, and inverse positivity gives $E^{n-1}\leq\varepsilon\Delta t$ whatever the policy. Letting $\varepsilon\downarrow0$ gives $E^{n-1}\leq0$, and backward induction proves the claim.

\subsection{Consistency of the discrete Hamiltonian and Bellman residual}\label{app:gfdm_consistency}
\begin{proposition}[Consistency of the discrete Hamiltonian and Bellman residual]\label{prop:consistency}
    Let $x_i$ carry a row of type (P), and let $R_i$ be the largest distance from $x_i$ to a point of its stencil. Then for $\varphi\in C^3$,
    \begin{equation*}
        \abs{\Delta_h\varphi(x_i)-\Delta\varphi(x_i)}\leq C\norm{\varphi}_{C^3}R_i.
    \end{equation*}
    If $b(x_i,m_i,\alpha)$ is represented by a retained fan sector, then for $\varphi\in C^2$,
    \begin{equation*}
        \abs*{(A_{\mathrm{fan}}^\alpha\pi_h\varphi)_i
        -b(x_i,m_i,\alpha)\cdot\nabla\varphi(x_i)}
        \leq CC_\theta\norm{\varphi}_{C^2}\abs{b(x_i,m_i,\alpha)}R_i.
    \end{equation*}
    Let $m_i$ be fixed, and suppose the discrete and continuum suprema below are attained in a common compact $K\subset\mcalA$. Then
    {\small
    \begin{equation*}
        \abs*{H_{i,h}[\pi_h\varphi;m_i]
        -H(x_i,\nabla\varphi(x_i),m_i)}
        \leq CC_\theta\norm{\varphi}_{C^2}R_i
        \sup_{\alpha\in K}\abs{b(x_i,m_i,\alpha)}.
    \end{equation*}}%
    Consequently, along a family with local support, if the preceding compact-control and drift bounds hold uniformly along the density histories $M$ of the family, then for test functions twice differentiable in time and three times in space,
    \begin{equation*}
        \abs{F_{h,\Delta t}^{n,i}(\pi_h\varphi;M)-F[\varphi;M_i^n](t_{n-1},x_i)}\leq C_\varphi(h+\Delta t)
    \end{equation*}
    holds uniformly over rows of type (P).
\end{proposition}

\begin{proof}
    Taylor expansion through degree two cancels under the exact Laplacian moments. Since the weights are nonnegative, the trace $\sum_{j\neq i}L_{ij}\abs{x_j-x_i}^2=4$ of the quadratic moments gives
    \begin{equation*}
        \sum_{j\neq i}L_{ij}\abs{x_j-x_i}^3\leq4R_i
    \end{equation*}
    for the cubic remainder. For the fan block, affine exactness cancels the constant and linear terms, while \autoref{lem:fan_existence} gives
    \begin{equation*}
        \sum_j\lambda_{ij}^\alpha\abs{z_{ij}}^2
        \leq R_i\sum_j\lambda_{ij}^\alpha\abs{z_{ij}}
        \leq C_\theta\abs{b(x_i,m_i,\alpha)}R_i.
    \end{equation*}
    The Hamiltonian estimate follows from
    \begin{equation*}
        \abs*{\sup_{\alpha\in K}f_\alpha-\sup_{\alpha\in K}g_\alpha}
        \leq\sup_{\alpha\in K}\abs{f_\alpha-g_\alpha}.
    \end{equation*}
    Implicit Euler contributes $O(\Delta t)$, the uniform compactness and drift bounds make the estimates uniform along $M$ under the Bellman supremum, and $R_i\leq\delta_i\leq C_\delta h$.
\end{proof}

\subsection{Consistency of mirror rows}\label{app:mirror_ghost_details}

Near a planar reflecting segment, choose tangential--normal coordinates with segment $\{r=0\}$, physical side $\{r\geq0\}$, and $\nu=-e_r$, so that $x_j=(y_j,s_{ij})$ and $x_j^{\mathrm{ghost}}=(y_j,-s_{ij})$ in the notation of \eqref{eq:bc_gamma_main}. For the homogeneous fold $u(x_j^{\mathrm{ghost}})=u(x_j)$, recall the folded weights $L_{ij}=\widetilde L_{ij}+\widetilde L^{\mathrm{ghost}}_{ij}$ of \autoref{sec:row_construction}, for each index $j$, the center included, and define the reflected physical-row operator by $(\Delta_h^R v)_i:=\sum_jL_{ij}v(x_j)$. The folded row acting on physical data equals the augmented row acting on the even extension.

We use one planar face and a one-sided physical pool. Without the obliqueness lower bound, positivity and $\msfP_2$ moments do not control the ghost-normal moment for a center lying $O(\delta_i)$ inside the boundary layer. On a finite cloud, minimizing over the feasible polytope without that bound can return $\widehat\gamma_{i,h}=0$. Because the added condition is linear in the weights, the mirror program remains a QP.

\begin{proposition}[Arbitrary-test expansion of a planar mirror row]\label{prop:mirror_boundary_expansion}
    Let $x_i$ carry a mirror row based on the plane of $\Gamma_0$, at distance at most $\kappa\delta_i$ from it, and let $x_0$ be its normal projection onto $\Gamma_0$. For all spatial $C^3$ test functions,
    \begin{equation}\label{eq:mirror_arbitrary_test}
        (\Delta_h^R\varphi)_i
        =\Delta\varphi(x_i)
        -\frac{\widehat\gamma_{i,h}}{\delta_i}
        \partial_\nu\varphi(x_0)+O(1)
    \end{equation}
holds uniformly along the refinement family.
    If the active controls lie in an $h$-independent compact set, then for $\varphi$ smooth in time and space the complete mirror-row Bellman residual satisfies
    \begin{equation}\label{eq:mirror_residual_expansion}
        F_{h,\Delta t}(\pi_h\varphi;M)(t_n,x_i)
        =F[\varphi;m](t_n,x_i)
        +\frac{\sigma^2}{2}\frac{\widehat\gamma_{i,h}}{\delta_i}
        \partial_\nu\varphi(t_n,x_0)+O(1).
    \end{equation}
\end{proposition}

\begin{proof}
    Choose a $C^3$ extension $\overline\varphi$ across the plane that equals $\varphi$ on the physical side. Every point of the augmented row lies within $(1+2\kappa)\delta_i$ of $x_i$. Augmented $\msfP_2$ exactness, nonnegativity and the trace $\sum_{k\neq i}\widetilde L_{ik}\abs{x_k-x_i}^2=4$ give
    \begin{equation*}
        \sum_k\widetilde L_{ik}\overline\varphi(x_k)
        =\Delta\varphi(x_i)+O(\delta_i).
    \end{equation*}
    Folding evaluates $\varphi$ at $x_j$ where the augmented row evaluates $\overline\varphi$ at $x_j^{\mathrm{ghost}}$, so $(\Delta_h^R\varphi)_i$ differs from this sum by $\sum_{j\in\mscrJ_i^{\mathrm{ghost}}}\widetilde L^{\mathrm{ghost}}_{ij}\bigl(\varphi(x_j)-\overline\varphi(x_j^{\mathrm{ghost}})\bigr)$. With $\nu=-e_r$, each difference is
    \begin{align*}
        \varphi(y_j,s_{ij})-\overline\varphi(y_j,-s_{ij})
        &=2s_{ij}\,\partial_r\varphi(0,0)
        +O(s_{ij}\abs{y_j}+s_{ij}^3)\\
        &=-2s_{ij}\partial_\nu\varphi(x_0)
        +O(s_{ij}\abs{y_j}+s_{ij}^3).
    \end{align*}
    The ghost of $x_j$ has displacement $(y_j,-s_{ij}-s_{ii})$ from $x_i$, so $s_{ij}\abs{y_j}+s_{ij}^3$ is at most a constant times the squared length of this displacement, and the trace bounds the remainder by $O(1)$. Summing the leading term gives exactly $-(\widehat\gamma_{i,h}/\delta_i)\partial_\nu\varphi(x_0)$, which completes the proof of \eqref{eq:mirror_arbitrary_test}.

    For a fixed bounded drift, compare the folded fan with the same positive ray representation applied to $\overline\varphi$. Each ghost fold contributes at most $C s_{ij}\lambda_{ij}^\alpha\norm{\varphi}_{C^1}$, and $s_{ij}\leq\abs{z_{ij}}$ with \autoref{lem:fan_existence} makes the total $O(\abs{b})$. These estimates are uniform over the compact active-control set and pass through the Bellman supremum. The implicit time difference contributes $O(\Delta t)$. Finally, the residual contains $-(\sigma^2/2)\Delta_h^R$, so the sign in \eqref{eq:mirror_residual_expansion} follows from \eqref{eq:mirror_arbitrary_test}.
\end{proof}

\subsection{Residual-to-value transfer and the smooth data rate}\label{app:smooth_rate}

For the pointwise-trace estimate below, write $\Lip_k(\overline\Omega)$ for the spatial Lipschitz class with constant $k$ independent of time.

For a continuous solution $u$, let $\pi_hu$ be its pointwise nodal sample and let $\pi_h^{\mathrm{cl}}u$ agree with it on dynamic and absorbing-Dirichlet nodes while using the copy rule \eqref{eq:copy_row} at copied nodes. Define
\begin{equation*}
    \epsilon_{\mathrm{res}}
    :=\max_{(n,i)\in\mscrI_h^{\mathrm{dyn}}}
    \abs*{F_{h,\Delta t}^{n,i}(\pi_h^{\mathrm{cl}}u;M)},
    \qquad
    \epsilon_{\mathrm{data}}
    :=\max_{(n,i)\in\mscrI_h^{\mathrm{data}}}
    \abs*{(\pi_h^{\mathrm{cl}}u)_i^n-g_{h,i}^n}.
\end{equation*}

\begin{corollary}[Residual-to-value comparison estimate]\label{cor:residual_value_transfer}
    Under the hypotheses of \autoref{thm:upwind_comparison}, let $U$ solve the discrete fan--Bellman problem and let $u\in C([0,T]\times\overline\Omega)$ be the viscosity solution of the limiting HJB problem for the prescribed density flow. Then
    \begin{equation*}
        \norm*{U-\pi_h^{\mathrm{cl}}u}_{\ell^\infty}
        \leq \epsilon_{\mathrm{data}}+T\epsilon_{\mathrm{res}}.
    \end{equation*}
    If, in addition, every copied node lies within $C_{\mathrm{cl}}h$ of its source node and $u(t,\cdot)\in\Lip_k(\overline\Omega)$ for all $t\in[0,T]$, then
    \begin{equation*}
        \norm*{U-\pi_hu}_{\ell^\infty}
        \leq \epsilon_{\mathrm{data}}+T\epsilon_{\mathrm{res}}
        +C_{\mathrm{cl}}kh.
    \end{equation*}
\end{corollary}

\begin{proof}
    Add and subtract
    $\epsilon_{\mathrm{data}}+\epsilon_{\mathrm{res}}(T-t_n)$
    from $\pi_h^{\mathrm{cl}}u$. The spatial generators annihilate constants, while the backward time difference shifts the residual by $\pm\epsilon_{\mathrm{res}}$. The resulting grid functions are respectively a discrete supersolution and subsolution, satisfy \eqref{eq:copy_row} with equality because the added constants are spatially uniform, and bracket the terminal and absorbing-Dirichlet data. \autoref{thm:upwind_comparison} gives
    \begin{equation*}
        \norm{U-\pi_h^{\mathrm{cl}}u}_{\ell^\infty}
        \leq\epsilon_{\mathrm{data}}+T\epsilon_{\mathrm{res}}.
    \end{equation*}
    The definition of $\Lip_k(\overline\Omega)$ gives $\norm{\pi_h^{\mathrm{cl}}u-\pi_hu}_{\ell^\infty}\leq C_{\mathrm{cl}}kh$, and the second estimate follows by the triangle inequality.
\end{proof}

\begin{remark}[First-order convergence on smooth data]\label{cor:smooth_rate}
    Assume \textup{(H2)}, let the density be prescribed and sampled nodally, with the Hamiltonian Lipschitz in time along it, let $u$ be twice continuously differentiable in time and three times in space on $[0,T]\times\overline\Omega$, let every dynamic row be of type (P) with a complete physical fan on a family with local support, let $\epsilon_{\mathrm{data}}\leq C_0(h+\Delta t)$, and let copied nodes lie within $C_{\mathrm{cl}}h$ of their sources and be read by no dynamic row. Then \autoref{prop:consistency} bounds the residual of $\pi_hu$ by $C(h+\Delta t)$, this is also the residual of $\pi_h^{\mathrm{cl}}u$, and \autoref{cor:residual_value_transfer} gives $\norm{U-\pi_hu}_{\ell^\infty}\leq C(h+\Delta t)$. The LQ benchmark takes exact boundary data and all its dynamic rows are of type (P), so on every cloud where the rows exist this holds with a constant that depends on the cloud only through the aperture cap and the ratio of the largest stencil radius to $h$.
\end{remark}

\subsection{Half-relaxed limits}\label{app:half_relaxed}

We write \textup{(C)} for comparison off the vertices: bounded upper semicontinuous subsolutions and bounded lower semicontinuous supersolutions of the terminal--boundary-value problem, with the terminal condition taken pointwise on $\overline\Omega\setminus(\Gamma_{\mathrm{c}}\cup\overline{\Gamma_{\mathrm{abs}}})$ and the lateral conditions in the relaxed viscosity sense on $\partial\Omega\setminus\Gamma_{\mathrm{c}}$ for $t<T$, are ordered on $Q_{\mathrm{c}}$.

\begin{lemma}[Half-relaxed limits at a prescribed density]\label{thm:viscosity_limit}\label{cor:frozen_density_convergence}\label{thm:planar_neumann_boundary}
    Let $d=2$, let the family have local support, and let $m\in C([0,T]\times\overline\Omega)$ be prescribed. Set $M_i^n:=m(t_n,x_i)$, and assume \textup{(H1)}--\textup{(H3)}. Let $\overline u$ and $\underline u$ be the upper and lower half-relaxed limits of the discrete solutions, which exist by \autoref{app:howard_solver}.
    \begin{enumerate}
        \item[(i)] On every $K\Subset\Omega$, $F_{h,\Delta t}(\pi_h\varphi;M)\to F[\varphi;m]$ for all smooth test functions, and $\overline u$ and $\underline u$ are respectively viscosity subsolutions and supersolutions in $[0,T)\times\Omega$.
        \item[(ii)] Let $\Gamma_0\Subset\Gamma_{\mathrm{refl}}$ be a relatively open planar segment at positive distance from $\Gamma_{\mathrm{abs}}$ and from all corners and boundary junctions, with outward unit normal $\nu$. Then, on $[0,T)\times\Gamma_0$, we have $\min\{F[\varphi;m],\partial_\nu\varphi\}\leq0$ at a local maximum of $\overline u-\varphi$ and $\max\{F[\varphi;m],\partial_\nu\varphi\}\geq0$ at a local minimum of $\underline u-\varphi$. Let $\Gamma_1\Subset\Gamma_{\mathrm{abs}}$ be relatively open and at positive distance from $\Gamma_{\mathrm{c}}$. Then, on $[0,T)\times\Gamma_1$, we have $\min\{F[\varphi;m],\overline u-g_{\mathrm{exit}}\}\leq0$ at such a maximum and $\max\{F[\varphi;m],\underline u-g_{\mathrm{exit}}\}\geq0$ at such a minimum.
        \item[(iii)] Assume \textup{(C)}. If $\overline u$ and $\underline u$ also satisfy the relaxed boundary inequalities at all points of $\partial\Omega\setminus\Gamma_{\mathrm{c}}$ not covered by \textup{(ii)}, and the terminal inequalities $\overline u(T,\cdot)\leq g_T\leq\underline u(T,\cdot)$ on $\overline\Omega\setminus(\Gamma_{\mathrm{c}}\cup\overline{\Gamma_{\mathrm{abs}}})$, then $\overline u=\underline u$ on $Q_{\mathrm{c}}$, where it coincides with the unique viscosity solution $u$ of the terminal--boundary-value problem, and $u_{h,\Delta t}\to u$ locally uniformly on $[0,T)\times(\overline\Omega\setminus\Gamma_{\mathrm{c}})$.
    \end{enumerate}
\end{lemma}

\subsection{Proof of \autoref{thm:viscosity_limit}\textup{(ii)} (planar reflecting segments and exit windows)}\label{app:proof_planar_neumann_boundary}

    Let $\overline u=\limsup{}^*u_{h,\Delta t}$, and suppose that $\overline u-\varphi$ has a strict local maximum at $(t_0,x_0)\in[0,T)\times\Gamma_0$. Since $t_0<T$ the index $n_h+1$ exists, and \eqref{eq:F_fan} reads the spatial level $n_h$ there, so the contact and the residual sit on the same spatial level. Data rows near $(t_0,x_0)$ are excluded by $t_0<T$ and by the reach of absorbing-Dirichlet rows in \autoref{def:row_classification}. Every copied value is attained at a dynamic row within $C_{\mathrm{cl}}h$, so the half-relaxed limit taken over all nodes coincides with the one taken over $\mscrI_h^{\mathrm{dyn}}$. Choose therefore a maximizing sequence $(t_{n_h},x_{i_h})\to(t_0,x_0)$ for $u_{h,\Delta t}-\pi_h\varphi$ \emph{over dynamic rows}. Copied nodes lie next to $\Gamma_{\mathrm{c}}$ and are out of reach, so, with the exclusion of data rows above, the entire stencil of the row at $x_{i_h}$ consists of dynamic rows for small $h$ and $\Delta t$. On it $u_{h,\Delta t}\leq\pi_h\varphi+c$ holds, with equality at the center, so no almost-maximality defect is incurred. Monotonicity, spatial constant invariance, and the discrete subsolution equation yield
    \begin{equation}\label{eq:boundary_contact_sub_discrete}
        F_{h,\Delta t}(\pi_h\varphi;M)(t_{n_h},x_{i_h})\leq o(1).
    \end{equation}

    Extract a subsequence of one row type. For type (P), \autoref{prop:consistency} applies. Its estimates use the exact nonnegative row, the fan and the support radius, and none of these asks the stencil ball to lie in $\Omega$. With the prescribed density continuous up to $\Gamma_0$, we obtain
    $F[\varphi;m](t_0,x_0)\leq0$. For type (R), the rest of $\partial\Omega$ lies at positive distance from $x_0$ and $\kappa\delta_{i_h}\to0$, so these rows reflect across the plane of $\Gamma_0$. Suppose instead that $\partial_\nu\varphi(t_0,x_0)>0$. Smoothness gives $\partial_\nu\varphi(t_{n_h},x_{0,h})\geq\eta>0$ at the normal projections $x_{0,h}$ for all sufficiently small $h$. By \autoref{prop:mirror_boundary_expansion}, non-degenerate diffusion and \eqref{eq:bc_gamma_main},
    \begin{equation*}
        F_{h,\Delta t}(\pi_h\varphi;M)(t_{n_h},x_{i_h})
        \geq \frac{\sigma^2}{2}\frac{\gamma_0\eta}{\delta_{i_h}}-C
        \longrightarrow+\infty,
    \end{equation*}
    contradicting \eqref{eq:boundary_contact_sub_discrete}. Hence a mirror subsequence forces $\partial_\nu\varphi(t_0,x_0)\leq0$. The two possible subsequences give
    \begin{equation*}
        \min\{F[\varphi;m],\partial_\nu\varphi\}(t_0,x_0)\leq0.
    \end{equation*}

    A general local maximum, interior or on $\Gamma_0$, reduces to the strict case by the standard perturbation $\varphi+\varepsilon(\abs{x-x_0}^2+(t-t_0)^2)$, which is even about the plane containing $\Gamma_0$ and so, unlike an odd tilt, adds nothing singular to \eqref{eq:mirror_residual_expansion}.

    For the lower half-relaxed limit, take a strict local minimum of $\underline u-\varphi$. The discrete supersolution equation reverses \eqref{eq:boundary_contact_sub_discrete}. A type (P) subsequence gives $F[\varphi;m](t_0,x_0)\geq0$ by the same application of \autoref{prop:consistency}, the stencil consisting of dynamic rows as above. On a type (R) subsequence, suppose $\partial_\nu\varphi(t_0,x_0)<0$. Smoothness gives $\partial_\nu\varphi(t_{n_h},x_{0,h})\leq-\eta<0$ at the normal projections for small $h$, which would make the singular term in \eqref{eq:mirror_residual_expansion} tend to $-\infty$, contradicting the supersolution inequality. Thus $\partial_\nu\varphi(t_0,x_0)\geq0$ on that branch, and
    \begin{equation*}
        \max\{F[\varphi;m],\partial_\nu\varphi\}(t_0,x_0)\geq0.
    \end{equation*}

    \emph{Exit windows.} Let $\overline u-\varphi$ have a strict local maximum at $(t_0,x_0)\in[0,T)\times\Gamma_1$, and take maxima $(t_{n_h},x_{i_h})\to(t_0,x_0)$ of $u_{h,\Delta t}-\pi_h\varphi$ over all nodes near $(t_0,x_0)$, with $u_{h,\Delta t}(t_{n_h},x_{i_h})\to\overline u(t_0,x_0)$. If they lie on data rows along a subsequence, \textup{(H1)} gives $\overline u(t_0,x_0)=g_{\mathrm{exit}}(t_0,x_0)$. Otherwise they lie on dynamic rows. For small $h$ the only face within reach of these nodes is the exit, so \autoref{def:row_classification} and \textup{(H3)} make these rows type (P). Their stencils lie among the nodes over which the maximum was taken, so monotonicity, spatial constant invariance and the discrete subsolution equation give $F_{h,\Delta t}(\pi_h\varphi;M)(t_{n_h},x_{i_h})\leq o(1)$, and \autoref{prop:consistency} gives $F[\varphi;m](t_0,x_0)\leq0$. This is the upper inequality on $\Gamma_1$. The perturbation above reduces a general local maximum to a strict one, and the same argument with the inequalities reversed gives the lower one.

\subsection{Attainment of the terminal data}\label{app:terminal_closure}

Parabolic comparison theorems with Neumann-type conditions assume the terminal inequality on $\overline\Omega$ \cite[Thm.~II.1]{Barles1993Neumann}. Some monotone treatments instead encode terminal data in the generalized boundary equation \cite{CalzolaCarliniDupuisSilva2023,JaroszkowskiJensen2024}. We prove the terminal inequality in $\Omega$ and on planar reflecting segments, which supplies the terminal hypothesis of \autoref{thm:viscosity_limit}\textup{(iii)} there.

\begin{lemma}[Attainment of the terminal data]\label{lem:terminal_closure}
    Let $m$, $M$, $\overline u$ and $\underline u$ be as in \autoref{thm:viscosity_limit}, and assume \textup{(H1)}--\textup{(H3)}. Then
    \begin{equation}\label{eq:terminal_closure}
        \overline u(T,x)\leq g_T(x)\leq\underline u(T,x)
    \end{equation}
    for $x\in\Omega$. If $\Gamma_0$ is as in \autoref{thm:viscosity_limit}\textup{(ii)}, then \eqref{eq:terminal_closure} also holds for $x\in\Gamma_0$.
\end{lemma}

\begin{proof}
    A time slope leaves a test function no contact below the terminal level, where the data are imposed. Let $C_U$ be the constant-barrier bound of \autoref{thm:viscosity_convergence}. Fix $x_0\in\Omega\cup\Gamma_0$, and choose $r\in(0,\min\{1,T\})$ so that $\overline B_{2r}(x_0)\Subset\Omega$ if $x_0\in\Omega$, and $\overline B_{2r}(x_0)\cap\partial\Omega\subset\Gamma_0$ if $x_0\in\Gamma_0$. Set $D:=\overline B_r(x_0)\cap\overline\Omega$, and for $\varepsilon>0$ let
    \begin{equation*}
        \varphi(t,x):=\frac{\abs{x-x_0}^2}{\varepsilon}+\theta\,(x-x_0)\cdot\nu+C(T-t),
    \end{equation*}
    with $\theta=0$ if $x_0\in\Omega$, $\theta=1$ if $x_0\in\Gamma_0$, and $C=C(\varepsilon)$ chosen below. Let $(\hat t,\hat x)$ maximize $\overline u-\varphi$ on $[T-r,T]\times D$. Comparing with the value at $(T,x_0)$ gives
    \begin{equation}\label{eq:terminal_contact}
        \overline u(T,x_0)\leq\overline u(\hat t,\hat x)-\frac{\abs{\hat x-x_0}^2}{\varepsilon}-\theta\,(\hat x-x_0)\cdot\nu-C(T-\hat t).
    \end{equation}
    Hence $\abs{\hat x-x_0}^2/\varepsilon+C(T-\hat t)\leq2C_U+1$, so for small $\varepsilon$ and large $C$ the maximum avoids the lateral and bottom faces of the cylinder. Add $\abs{x-\hat x}^4+(t-\hat t)^2$ to $\varphi$ to make the maximum strict. This changes neither its time derivative nor its spatial derivatives up to order two at $(\hat t,\hat x)$, nor, when $\hat x\in\Gamma_0$, its normal derivative on the plane of $\Gamma_0$. A copied value is attained at a dynamic or terminal row within $C_{\mathrm{cl}}h$, so along a subsequence realizing $\overline u(\hat t,\hat x)$ the maxima $(t_{n_h},x_{i_h})$ of $u_{h,\Delta t}-\pi_h\varphi$ over the dynamic and terminal rows in $[T-r,T]\times D$ converge to $(\hat t,\hat x)$, and $u_{h,\Delta t}(t_{n_h},x_{i_h})\to\overline u(\hat t,\hat x)$.

    Suppose $t_{n_h}<T$ along a further subsequence. Absorbing-Dirichlet rows are out of reach of $\hat x$. Copied nodes lie next to $\Gamma_{\mathrm{c}}$, so they are out of reach too. So the residual at $(t_{n_h},x_{i_h})$ reads dynamic rows at level $t_{n_h}$ and, at level $t_{n_h+1}$, only its center, a dynamic or terminal row. All of them lie in the set over which the maximum was taken, so, as in \autoref{app:proof_planar_neumann_boundary}, monotonicity, spatial constant invariance and the discrete subsolution equation give $F_{h,\Delta t}(\pi_h\varphi;M)(t_{n_h},x_{i_h})\leq o(1)$. On type (P) rows, \autoref{prop:consistency} gives $F[\varphi;m](\hat t,\hat x)\leq0$, that is,
    \begin{equation*}
        C-\frac{\sigma^2d}{\varepsilon}+H\parens[\big]{\hat x,\nabla\varphi(\hat t,\hat x),m(\hat t,\hat x)}\leq0,
    \end{equation*}
    which fails once $C$ exceeds $\sigma^2d/\varepsilon$ plus the bound of $\abs{H}$ over $D$, the range of $m$ and gradients up to $2r/\varepsilon+2$, finite by \textup{(H2)}. Type (R) rows occur only when $\hat x\in\Gamma_0$, since no face is within reach near $\hat x\in\Omega$. There $\partial_\nu\varphi=1$ on the plane of $\Gamma_0$, so \eqref{eq:mirror_residual_expansion} and \eqref{eq:bc_gamma_main} send the residual to $+\infty$, a contradiction. Hence $t_{n_h}=T$ along the subsequence, so $\hat t=T$, and \textup{(H1)} gives $\overline u(T,\hat x)=\lim g_{h,i_h}^{N_t}=g_T(\hat x)$. Then \eqref{eq:terminal_contact} gives $\overline u(T,x_0)\leq g_T(\hat x)+\theta\abs{\hat x-x_0}$, and letting $\varepsilon\downarrow0$, so that $\hat x\to x_0$, proves the upper inequality. The lower one follows from the same argument with $-\varphi$, for which the discrete supersolution equation reverses the inequalities and the singular term tends to $-\infty$.
\end{proof}

When $g_{h,i}^{N_t}=G(x_i,M^{N_t})$, with $M^{N_t}$ the nodal sample of the prescribed $m(T)$, the consistency in \textup{(H1)} follows from the continuity of $G$.

\subsection{Proof of \autoref{thm:viscosity_convergence} and of \autoref{thm:viscosity_limit}\textup{(i)} and \textup{(iii)}}\label{app:proof_viscosity_convergence}

    \emph{Monotonicity.} At each dynamic row \textup{(H3)} gives nonnegative off-center Laplacian and fan coefficients with zero row sum. Hence for any policy the fixed-policy matrix is an M-matrix by \autoref{thm:constructed_operator_monotone}(ii), and the supremum of the corresponding affine residuals is monotone in the Barles--Souganidis sense.

    \emph{Consistency.} Let $K\Subset\Omega$. For small $h$ no face is within reach of the nodes approaching $K$, so \textup{(H3)} gives them rows of type (P). \autoref{prop:consistency} gives convergence of the diffusion and fan blocks on smooth test functions. Compactness or coercivity confines the active controls on a bounded test gradient, and the estimate therefore passes uniformly through the Bellman supremum. Since $M_i^n=m(t_n,x_i)$ with $m$ uniformly continuous on $[0,T]\times\overline\Omega$, and $b,\mcalL$ are uniformly continuous on the compact active set of \textup{(H2)}, $H(x_i,p,M_i^n)\to H(x_0,p,m(t_0,x_0))$ along the grid points, uniformly for bounded $p$, which identifies the continuous Hamiltonian.

    \emph{Stability.} Let $C_0$ bound the imposed terminal and absorbing-Dirichlet data. Because both spatial generators annihilate constants on interior and homogeneous reflecting rows, the grid functions
    \begin{equation*}
        \overline U_i^n=C_0+C_H(T-t_n),
        \qquad
        \underline U_i^n=-C_0-C_H(T-t_n)
    \end{equation*}
    are respectively a discrete supersolution and subsolution. Indeed, on a constant state the row residual is the time slope plus $\sup_{\alpha\in\mcalA}\{-\mcalL(x_i,M_i^n,\alpha)\}$, whose absolute value is bounded by $C_H$. The discrete comparison principle is \autoref{thm:upwind_comparison}, whose hypotheses \textup{(H3)} provides, so it yields $\abs{U_i^n}\leq C_0+C_HT$, uniformly in $h$ and $\Delta t$.

    \emph{Half-relaxed limits.} This uniform bound makes
    $\overline u=\limsup{}^*u_h$ and $\underline u=\liminf{}_*u_h$ finite. If $\overline u-\varphi$ has a strict local maximum at $(t_0,x_0)\in[0,T)\times\Omega$, choose nearby discrete maxima of $u_h-\pi_h\varphi$. Monotonicity of the same residual \eqref{eq:F_fan}, constant invariance of its spatial part, and consistency give $F[\varphi;m](t_0,x_0)\leq0$. The symmetric minimum argument gives the supersolution inequality for $\underline u$. Thus these are interior viscosity subsolutions and supersolutions.

    \autoref{lem:terminal_closure} gives the terminal inequalities in $\Omega$. Under the hypotheses of \autoref{thm:viscosity_limit}\textup{(iii)}, comparison \textup{(C)} gives $\overline u\leq\underline u$ on $Q_{\mathrm{c}}$. The reverse inequality is automatic. Hence the limits coincide and the discrete solutions converge locally uniformly. Without spatial boundary identification the argument yields only the interior subsolution and supersolution conclusions.

\subsection{Comparison off the boundary vertices}\label{app:comparison}

Neumann-type comparison theorems assume $\partial\Omega$ of class $W^{3,\infty}$ \cite[Thm.~II.1]{Barles1993Neumann}, or of class $C^1$ under stronger conditions on $H$ \cite{Ishii1991Oblique}. The nonsmooth result \cite[Thm.~2.1]{DupuisIshii1990} is elliptic, under strict monotonicity in $u$ and a uniform exterior cone along the reflection direction, which the outward normal does not provide at a re-entrant corner. Among monotone schemes, the finite element treatment of mixed conditions assumes comparison outright \cite[Assump.~8.2]{JaroszkowskiJensen2024}, and the game-based scheme needs a $C^2$ domain \cite{Daniel2013Game}. None of these provides \textup{(C)} on the evacuation geometry. Convexity in the gradient does, because it lets a strict subsolution that tends to $-\infty$ remove the boundary vertices and the terminal exit from the comparison.

\begin{proposition}[Comparison for convex Hamiltonians of at most quadratic growth]\label{prop:comparison}
    Let $\Omega\subset\mbbR^2$ be bounded and lie locally on one side of its boundary, and let $\partial\Omega\setminus\Gamma_{\mathrm{c}}$ be of class $C^2$. Near each point of $\Gamma_{\mathrm{c}}$ let $\partial\Omega$ consist of two straight segments ending at that point, at an interior angle in $(0,2\pi)$, let $\Gamma_{\mathrm{abs}}\setminus\Gamma_{\mathrm{c}}$ be a finite union of open straight segments, and let each point of $\Gamma_{\mathrm{refl}}\cap\overline{\Gamma_{\mathrm{abs}}}$ lie in the relative interior of a straight segment of $\partial\Omega$. Let $g_T$ and $g_{\mathrm{exit}}$ be bounded, with $g_{\mathrm{exit}}$ continuous on $[0,T)\times\Gamma_{\mathrm{abs}}$, and assume \textup{(H4)}. Then \textup{(C)} holds.
\end{proposition}

Every Hamiltonian of the form \eqref{eq:bellman_hamiltonian} is convex in $p$, being a supremum of affine functions, and convexity with \eqref{eq:comparison_structure} makes $H_m$ continuous. The other conditions hold in two cases.
\begin{itemize}
    \item \emph{Separable.} $H_m(t,x,p)=H_0(p)-f(t,x)$, with $H_0$ convex of at most quadratic growth and $f$ continuous. Then $L=0$ and $\omega_H$ is the modulus of $f$.
    \item \emph{Compact controls.} $\mcalA$ is compact, $b(x,m(t,x),\alpha)$ is bounded and Lipschitz in $(t,x)$, and $\mcalL(x,m(t,x),\alpha)$ is uniformly continuous in $(t,x)$, each uniformly in $\alpha$. Then $H_m$ grows linearly.
\end{itemize}
The growth condition \eqref{eq:comparison_growth} is an addition to \textup{(H2)}. The evacuation game is separable, with $H_0=\abs{p}^2/2$ and $f=V_{\mathrm{wall}}+V_{\mathrm{obs}}+\zeta m$ for a prescribed continuous $m$. Its domain meets the geometric hypotheses: the box corners and wall feet have interior angle $\pi/2$, the door jambs $3\pi/2$ and the exit ends $\pi$, the exit windows lie inside straight outer-wall segments, and the circle is a $C^2$ curve away from $\Gamma_{\mathrm{c}}$. The LQ benchmark is also separable at a prescribed density, and its square frame, absorbing on every side, has no junction of reflecting and absorbing parts.

\begin{proof}
    Let $u$ be a subsolution and $v$ a supersolution, and let $\Sigma_{\mathrm{c}}:=([0,T]\times\Gamma_{\mathrm{c}})\cup(\{T\}\times\overline{\Gamma_{\mathrm{abs}}})$ be the complement of $Q_{\mathrm{c}}$ in $[0,T]\times\overline\Omega$. We build a strict subsolution $\varphi$ that tends to $-\infty$ at $\Sigma_{\mathrm{c}}$, so that $\theta u-v+(1-\theta)\varphi$ attains its maximum off $\Sigma_{\mathrm{c}}$, and we exclude a positive maximum by doubling variables. The limit $\theta\uparrow1$ then gives $u\leq v$ on $Q_{\mathrm{c}}$.

    \emph{Step 1: a logarithmic identity.} A supporting plane of $H_m$ at $p=0$ has slope at most $2\sup_{\abs{p}\leq1}\abs{H_m}$, which is finite. Enlarging $C_{\mathrm{q}}$, we may therefore also assume $H_m\geq-C_{\mathrm{q}}(1+\abs{p})$. Put $\varkappa:=\sigma^2/(2C_{\mathrm{q}})$. For a smooth positive $\psi$ and $\varphi:=-\varkappa\log\psi$,
    \begin{equation*}
        -\partial_t\varphi-\frac{\sigma^2}{2}\Delta\varphi+C_{\mathrm{q}}\abs{\nabla\varphi}^2=\frac{\varkappa}{\psi}\parens[\Big]{\partial_t\psi+\frac{\sigma^2}{2}\Delta\psi},
        \qquad
        \partial_\nu\varphi=-\varkappa\frac{\partial_\nu\psi}{\psi},
    \end{equation*}
    so \eqref{eq:comparison_growth} gives
    \begin{equation}\label{eq:comparison_log}
        F[\varphi;m]\leq\frac{\varkappa}{\psi}\parens[\Big]{\partial_t\psi+\frac{\sigma^2}{2}\Delta\psi}+C_{\mathrm{q}}.
    \end{equation}

    \emph{Step 2: the Dirichlet data are attained.} Let $x_0\in\Gamma_{\mathrm{abs}}\setminus\Gamma_{\mathrm{c}}$ and $t_0<T$. Let $r$ be the inward distance to the line carrying $x_0$, fix $A$ with $\sigma^2A>C_{\mathrm{q}}$, and put
    \begin{equation*}
        \Theta:=K(r-Ar^2)+K_2\parens[\big]{\abs{x-x_0}^4+(t-t_0)^2}
    \end{equation*}
    on a cylinder of radius $\varrho\leq1/(4A)$ about $(t_0,x_0)$, with $t_0+\varrho<T$, $\varrho<\operatorname{dist}(x_0,\Gamma_{\mathrm{c}})$ and $\overline{B(x_0,\varrho)}\cap\partial\Omega\subset\Gamma_{\mathrm{abs}}$. There $r-Ar^2\geq0$ and $\abs{\nabla(r-Ar^2)}\leq1$, and $-\frac{\sigma^2}{2}\Delta\Theta\mp\partial_t\Theta\geq\sigma^2AK-C(K_2,\varrho)$, with $C(K_2,\varrho)$ independent of $K$. A contact at $t=0$ is covered by the solution concept.

    Suppose $u(t_0,x_0)>g_{\mathrm{exit}}(t_0,x_0)$. Choose $\varrho$ so that $g_{\mathrm{exit}}$ varies by less than the difference, then $K_2$, independent of $K$ because $K(r-Ar^2)\geq0$, so that $u-\Theta$ cannot attain its maximum where $\abs{x-x_0}=\varrho$ or $\abs{t-t_0}=\varrho$. At the maximum point $u\geq u(t_0,x_0)$, so $u>g_{\mathrm{exit}}$, and the relaxed condition gives $F[\Theta;m]\leq0$. The lower bound on $H_m$ gives $F[\Theta;m]\geq K(\sigma^2A-C_{\mathrm{q}})-C(K_2,\varrho)$, which is positive for $K$ large, a contradiction. Hence $u\leq g_{\mathrm{exit}}$ on $[0,T)\times(\Gamma_{\mathrm{abs}}\setminus\Gamma_{\mathrm{c}})$.

    Suppose $v(t_0,x_0)<g_{\mathrm{exit}}(t_0,x_0)$. Put $w:=e^{-v/\varkappa}$, which is upper semicontinuous with $w\leq M_v:=e^{\norm{v}_\infty/\varkappa}$, and $\Theta_0:=w(t_0,x_0)+\Theta$. Choose $\varrho$ as before and $K_2$, again independent of $K$, so that $w-\Theta_0$ attains its maximum $c'\geq0$ away from those faces. There $w\geq\Theta_0\geq w(t_0,x_0)$, so $v\leq v(t_0,x_0)<g_{\mathrm{exit}}$. The function $\xi:=-\varkappa\log(\Theta_0+c')$ touches $v$ from below at that point, so the relaxed condition gives $F[\xi;m]\geq0$. For $K$ large the bracket in \eqref{eq:comparison_log} is negative, and $\Theta_0+c'=w\leq M_v$ at the contact, so
    \begin{equation*}
        F[\xi;m]\leq-\frac{\varkappa}{M_v}\parens[\big]{\sigma^2AK-C(K_2,\varrho)}+C_{\mathrm{q}},
    \end{equation*}
    which is negative for $K$ large, a contradiction. Hence $v\geq g_{\mathrm{exit}}$ on $[0,T)\times(\Gamma_{\mathrm{abs}}\setminus\Gamma_{\mathrm{c}})$.

    \emph{Step 3: the barrier.} Choose $r_0>0$ so small that the balls $B(x_c,2r_0)$, $x_c\in\Gamma_{\mathrm{c}}$, are disjoint and each meets $\partial\Omega$ only in the two straight segments at $x_c$. Let $\beta\in C^\infty((0,\infty))$ be nonnegative, with $\beta(r)=\log(r_0/r)$ on $(0,r_0/2)$ and $\beta=0$ on $[r_0,\infty)$, and set $\Psi_{\mathrm{c}}(x):=\sum_{x_c\in\Gamma_{\mathrm{c}}}\beta(\abs{x-x_c})$. Then $\Psi_{\mathrm{c}}\geq0$ tends to $+\infty$ at $\Gamma_{\mathrm{c}}$, $\Delta\Psi_{\mathrm{c}}$ is bounded because $\log r$ is harmonic in two dimensions, and $\partial_\nu\Psi_{\mathrm{c}}=0$ on $\partial\Omega\setminus\Gamma_{\mathrm{c}}$, because a radial function has zero normal derivative on a straight segment through its center.

    Let $\Gamma$ be a connected component of $\Gamma_{\mathrm{abs}}\setminus\Gamma_{\mathrm{c}}$, lying on a line $\ell$, and let $\varepsilon>0$ be so small that the $3\varepsilon$-neighborhood of $\overline\Gamma$ meets $\Gamma_{\mathrm{refl}}$ only in $\ell$. Let $\chi_1\geq0$ be smooth, with $\chi_1=1$ near $0$ and $\chi_1=0$ on $[\varepsilon,\infty)$, and put
    \begin{equation*}
        \beta_\Gamma(x):=\int_\Gamma\frac{\chi_1(\abs{x-z})}{\abs{x-z}}\diff z,
        \qquad
        \vartheta_\Gamma(\tau,x):=\mbbE\bracks[\big]{\beta_\Gamma(x+\sigma\sqrt\tau Z)},
        \qquad x\in\mbbR^2,
    \end{equation*}
    with $Z$ standard normal in $\mbbR^2$. The function $\beta_\Gamma\geq0$ is integrable, smooth off $\overline\Gamma$, even across $\ell$, supported in the $\varepsilon$-neighborhood of $\overline\Gamma$, and tends to $+\infty$ at $\overline\Gamma$. Hence $\vartheta_\Gamma\geq0$ solves $\partial_\tau\vartheta_\Gamma=\frac{\sigma^2}{2}\Delta\vartheta_\Gamma$, is even across $\ell$, is smooth up to $\tau=0$ off $\overline\Gamma$, and tends to $+\infty$ as $(\tau,x)\to(0,z)$, $z\in\overline\Gamma$, by Fatou's lemma. Let $\chi\in C_c^\infty(\mbbR^2)$ be nonnegative and even across $\ell$, with $\chi=1$ on the $2\varepsilon$-neighborhood of $\overline\Gamma$ and with support in the $3\varepsilon$-neighborhood, and set $\Psi_\Gamma(t,x):=\chi(x)\vartheta_\Gamma(T-t,x)$. Then $\partial_\nu\Psi_\Gamma=0$ on $\Gamma_{\mathrm{refl}}\setminus\overline\Gamma$, and
    \begin{equation*}
        -\partial_t\Psi_\Gamma-\frac{\sigma^2}{2}\Delta\Psi_\Gamma=-\frac{\sigma^2}{2}\parens[\big]{2\nabla\chi\cdot\nabla\vartheta_\Gamma+\vartheta_\Gamma\Delta\chi},
    \end{equation*}
    which is bounded, because $\nabla\chi$ is supported at distance at least $\varepsilon$ from the support of $\beta_\Gamma$, where the Gaussian kernel and its gradient are bounded uniformly in $\tau$.

    Set $\Psi:=1+\Psi_{\mathrm{c}}+\sum_\Gamma\Psi_\Gamma$, $\psi:=e^{\Lambda(T-t)}\Psi$, $\varphi:=-\varkappa\log\psi-C_0$ and $\lambda_0:=(1+C_{\mathrm{q}})/\varkappa$. Since $\Psi\geq1$ and $\partial_t\Psi+\frac{\sigma^2}{2}\Delta\Psi$ is bounded,
    \begin{equation*}
        -\partial_t\psi-\frac{\sigma^2}{2}\Delta\psi-\lambda_0\psi=e^{\Lambda(T-t)}\parens[\Big]{(\Lambda-\lambda_0)\Psi-\partial_t\Psi-\frac{\sigma^2}{2}\Delta\Psi}\geq0
    \end{equation*}
    for $\Lambda$ large, and \eqref{eq:comparison_log} gives $F[\varphi;m]\leq-1$ off $\Sigma_{\mathrm{c}}$. Moreover, $\varphi$ is smooth off $\Sigma_{\mathrm{c}}$, $\varphi\to-\infty$ at $\Sigma_{\mathrm{c}}$, and $\partial_\nu\varphi=0$ on $\Gamma_{\mathrm{refl}}\setminus\Gamma_{\mathrm{c}}$. Choose $C_0$ so that $\varphi\leq g_T$ and $\varphi\leq g_{\mathrm{exit}}$.

    \emph{Step 4: the convex combination.} Suppose $u>v$ at a point of $Q_{\mathrm{c}}$. Fix $\theta\in(0,1)$ so close to $1$ that $\Phi_0:=\theta u-v+(1-\theta)\varphi$ is positive there. Extended by $-\infty$ on $\Sigma_{\mathrm{c}}$, this function is upper semicontinuous on $[0,T]\times\overline\Omega$, so it attains a positive maximum at some $(\hat t,\hat x)\in Q_{\mathrm{c}}$. If $\hat t=T$, the terminal condition gives $\Phi_0\leq(1-\theta)(\varphi-g_T)\leq0$. If $\hat x\in\Gamma_{\mathrm{abs}}$, Step 2 gives $\Phi_0\leq(1-\theta)(\varphi-g_{\mathrm{exit}})\leq0$. So $\hat t<T$, and $\hat x$ lies in $\Omega$ or on $\Gamma_{\mathrm{refl}}\setminus\Gamma_{\mathrm{c}}$.

    Put $\varpi(t,x):=\abs{x-\hat x}^4+(t-\hat t)^2$. If $\hat x\in\Gamma_{\mathrm{refl}}$, let $n:=\nu(\hat x)$ and $\rho>0$. If $\hat x\in\Omega$, take $\rho=0$. On a small closed neighborhood of $(\hat t,\hat x)$ in $[0,T)\times\overline\Omega$, maximize
    \begin{equation*}
        \Phi:=\theta u(s,y)-v(t,x)+(1-\theta)\varphi(t,x)-\varpi(t,x)-\frac{\abs{x-y}^2+(t-s)^2}{2\eta}-\rho\,n\cdot(x+y-2\hat x).
    \end{equation*}
    As $\eta\to0$, maximizers converge along a subsequence to a maximizer $(\bar t,\bar x)$ of $\Phi_0-\varpi-2\rho\,n\cdot(x-\hat x)$, with $(\abs{x-y}^2+(t-s)^2)/\eta\to0$ \cite[Lem.~3.1]{CrandallIshiiLions1992}. As $\rho\to0$, $(\bar t,\bar x)\to(\hat t,\hat x)$, and since $\Phi_0$ is largest at $(\hat t,\hat x)$, we have $\varpi(\bar t,\bar x)\leq2\rho\abs{n\cdot(\bar x-\hat x)}$.

    The test gradients are $q_u:=((y-x)/\eta+\rho n)/\theta$ for $u$ at $(s,y)$, and
    \begin{equation*}
        q_v=\theta q_u+(1-\theta)\nabla\widehat\varphi(t,x),
        \qquad
        \widehat\varphi:=\varphi-\frac{\varpi+2\rho\,n\cdot x}{1-\theta},
    \end{equation*}
    for $v$ at $(t,x)$. The normal $n$ is frozen, so that $x$ and $y$ see opposite gradients of the penalty and the remainder of $q_v$ falls on $\widehat\varphi$, whose gradient stays bounded.

    Now take $x$ or $y$ on $\Gamma_{\mathrm{refl}}$, within a neighborhood where $n\cdot\nu\geq\frac12$. The exterior sphere property of a $C^2$ boundary \cite[(7.13), (7.21)]{CrandallIshiiLions1992} gives $(x-y)\cdot\nu(y)\leq C\abs{x-y}^2$ and $(y-x)\cdot\nu(x)\leq C\abs{x-y}^2$, so we have $\theta\,\nu(y)\cdot q_u\geq\rho\,n\cdot\nu(y)-C\abs{x-y}^2/\eta>0$, and, using $\partial_\nu\varphi=0$,
    \begin{equation*}
        \nu(x)\cdot q_v\leq\abs{\nabla\varpi(t,x)}+\frac{C\abs{x-y}^2}{\eta}-\rho\,n\cdot\nu(x)<0
    \end{equation*}
    for small $\eta$ and $\rho$. Here $\bar x$ on $\Gamma_{\mathrm{refl}}$ has $\abs{n\cdot(\bar x-\hat x)}\leq C\abs{\bar x-\hat x}^2$, so the bound on $\varpi$ gives $\abs{\bar x-\hat x}^2\leq C\rho$ and $\abs{\nabla\varpi}\leq C\rho^{3/2}$. If $\bar x\in\Omega$, then $x$ and $y$ lie in $\Omega$ for small $\eta$. By the continuity of $\nu$, these signs persist at the approximating points of the closed jets, so the relaxed conditions yield the equations on both sides.

    Apply \cite[Thm.~3.2]{CrandallIshiiLions1992} on the locally compact set $[0,T)\times\overline\Omega$ to $\theta u-\rho\,n\cdot y$ and $v-(1-\theta)\varphi+\varpi+\rho\,n\cdot x$. This gives $X\leq Y$ for the two spatial Hessians, with the same time slope $(s-t)/\eta$, and it covers $\hat t=0$ under the convention of the solution concept. Subtract $\theta$ times the subsolution inequality from the supersolution inequality. The time slope and the Hessians combine to the derivatives of $(1-\theta)\widehat\varphi$, and the term $\frac{\sigma^2}{2}\operatorname{tr}(X-Y)\leq0$ is dropped. The convexity of $H_m$ gives $H_m(t,x,q_v)\leq\theta H_m(t,x,q_u)+(1-\theta)H_m(t,x,\nabla\widehat\varphi)$. Therefore
    \begin{equation*}
        0\leq(1-\theta)F[\widehat\varphi;m](t,x)+\theta\parens[\big]{H_m(t,x,q_u)-H_m(s,y,q_u)}.
    \end{equation*}
    At fixed $\theta$, the derivatives of $(\varpi+2\rho\,n\cdot x)/(1-\theta)$ at $(\bar t,\bar x)$ tend to $0$ as $\rho\to0$, so $F[\widehat\varphi;m]\leq-\frac12$ near $(\bar t,\bar x)$ for small $\rho$. By \eqref{eq:comparison_structure} and $\abs{q_u}\leq(\abs{x-y}/\eta+\rho)/\theta$, the last term tends to $0$ as $\eta\to0$, because $(\abs{t-s}+\abs{x-y})\abs{x-y}/\eta\leq2(\abs{t-s}^2+\abs{x-y}^2)/\eta$. This gives $0\leq-(1-\theta)/2$, a contradiction.
\end{proof}

The convex combination is that of \cite{DaLioLey2006}, who work on $\mbbR^N$. The frozen normal and the quartic localization follow the proof of \cite[Thm.~7.5]{CrandallIshiiLions1992}. The barrier at the vertices is the PDE form of \cite[Thm.~2.2]{VaradhanWilliams1985}. For normal reflection, their parameter $\alpha$ vanishes on a wedge of any angle in $(0,2\pi)$, and Brownian motion started away from the vertex does not reach it. Boundedness cannot be dropped. On a domain star-shaped about a vertex $x_c$, with all sides reflecting, $g_T\equiv0$ and $H_m=\abs{p}^2/(2c)$, the function $u\equiv0$ is a solution, and, with $E_1(z):=\int_z^\infty e^{-s}s^{-1}\diff s$, $\varkappa=c\sigma^2$ and any $a>0$,
\begin{equation*}
    v(t,x)=-\varkappa\log\parens[\Big]{1+aE_1\parens[\Big]{\frac{\abs{x-x_c}^2}{2\sigma^2(T-t)}}}
\end{equation*}
is a supersolution with the same terminal value off $x_c$, unbounded below near $x_c$ for $t<T$, and $v<u$.

\subsection{Proof of \autoref{thm:polygon}}\label{app:polygon}

Call a closed straight segment of $\partial\Omega$ with endpoints in $\Gamma_{\mathrm{c}}$ and no point of $\Gamma_{\mathrm{c}}$ in its relative interior a \emph{side}, and write $\alpha_v$ for the interior angle at $v\in\Gamma_{\mathrm{c}}$.

\begin{lemma}[Balls away from the vertices]\label{lem:polygon_balls}
    Under the geometric hypotheses of \autoref{thm:polygon} there are $C_\Gamma\geq1$ and $r_0>0$ such that, for $r<r_0$ and $x\in\overline\Omega$ with $\operatorname{dist}(x,\Gamma_{\mathrm{c}})\geq C_\Gamma r$, the intersection $B(x,r)\cap\partial\Omega$ is empty or is a chord of the ball on one side.
\end{lemma}

\begin{proof}
    Let $2r_0$ be the least distance between two disjoint sides, with $r_0:=1$ if there are none. Let $y$ and $z$ lie on two sides with common endpoint $v$, put $\beta_v:=\min(\alpha_v,2\pi-\alpha_v)$, and put $\sigma_v:=\sin\beta_v$ if $\beta_v<\pi/2$ and $\sigma_v:=1$ otherwise. The distance from $y$ to the second side is at least $\sigma_v\abs{y-v}$, so $\abs{y-z}\geq\sigma_v\abs{y-v}$. If $B(x,r)$, $r<r_0$, met two sides at $y$ and $z$, then $\abs{y-z}<2r<2r_0$ would make them adjacent, and then $\abs{x-v}\leq r+\abs{y-v}<r+2r/\sigma_v$. So $C_\Gamma:=1+2/\min_v\sigma_v$ excludes two sides. Since $C_\Gamma\geq1$, the ball contains no endpoint, and its intersection with the one side it meets is a chord.
\end{proof}

Put $C_{\mathrm{v}}:=C_\Gamma C_\delta$.

\begin{lemma}[Copy sources]\label{lem:copy_source}
    Under the hypotheses of \autoref{thm:polygon} there are $C_{\mathrm{cl}},h_0>0$ such that, for $h<h_0$, each interior node within $C_{\mathrm{v}}h$ of $\Gamma_{\mathrm{c}}$ to which none of (P), (D), (R) applies has a node with a dynamic row within $C_{\mathrm{cl}}h$.
\end{lemma}

\begin{proof}
    Let $v\in\Gamma_{\mathrm{c}}$. If both sides at $v$ are absorbing, a node near $v$ at distance at least $C_rh_i$ from $\partial\Omega$ receives a row of type (P), and a nearer node has an absorbing nearest face within $\kappa\delta_i$, so it receives a row of type (P) or (D). Otherwise let $e$ be the unit vector bisecting the interior angle at $v$. Near $v$ the domain is a wedge $W$ with apex $v$. If $x\in W$, then $x+te\in W$ for $t\geq0$. For $\alpha_v\leq\pi$ this holds because $W$ is a convex cone containing $e$. For $\alpha_v>\pi$, the complement of $W$ is a closed convex cone $K$ with apex $v$ containing the direction $-e$, so $x+te\in K$ would give $x\in K$.

    Let $x_i$ be near $v$. The point $y:=x_i+\frac12\delta_ie$ lies in $B(x_i,\delta_i)\cap\Omega$ and is visible from $x_i$, so \autoref{def:stencil_geometry} gives $x_j\in\mcalV_i$ with $\abs{x_j-y}\leq h$. Hence
    \begin{equation*}
        (x_j-v)\cdot e\geq(x_i-v)\cdot e+\parens[\Big]{\frac{c_\delta}{2}-1}h,
        \qquad
        \abs{x_j-x_i}\leq C_\delta h.
    \end{equation*}
    The node $x_j$ is interior, because no reflecting node is sampled and $y$ lies farther than $h$ from any exit, so the step repeats. From $(x_i-v)\cdot e\geq-C_{\mathrm{v}}h$, a number of steps depending only on $c_\delta,C_{\mathrm{v}},C_r$ reaches a node with $(x-v)\cdot e\geq(C_{\mathrm{v}}+C_r)h$. For small $h$ this walk of bounded length stays in $W$ and away from the other points of $\Gamma_{\mathrm{c}}$, so its end lies outside the $C_{\mathrm{v}}h$-neighborhood of $\Gamma_{\mathrm{c}}$, and \autoref{lem:polygon_balls} and \autoref{lem:mirror_constructibility} give it a row. If both sides at $v$ are reflecting, this row is dynamic. Otherwise one side is absorbing and $\alpha_v=\pi$ by the hypothesis on $\Gamma_{\mathrm{refl}}\cap\overline{\Gamma_{\mathrm{abs}}}$, so $(x-v)\cdot e$ is the distance to the line of the sides at $v$, $B(x,C_rh_i)\subset\Omega$, and the row is of type (P).
\end{proof}

    We now prove \autoref{thm:polygon}. Take $h$ so small that $C_\delta h<r_0$. A node with $\operatorname{dist}(x_i,\Gamma_{\mathrm{c}})\geq C_{\mathrm{v}}h$ has $\operatorname{dist}(x_i,\Gamma_{\mathrm{c}})\geq C_\Gamma\delta_i$, so \autoref{lem:polygon_balls} and \autoref{lem:mirror_constructibility} give it a row of type (P), (D) or (R). A nearer node takes such a row or the copy row that \textup{(H3)} allows, whose source exists by \autoref{lem:copy_source}. So \textup{(H3)} holds.

    A point of $\partial\Omega\setminus\Gamma_{\mathrm{c}}$ lies in the relative interior of a side, so it lies in a segment $\Gamma_0$ or a set $\Gamma_1$ of \autoref{thm:viscosity_limit}\textup{(ii)}. So the relaxed boundary inequalities hold on all of $[0,T)\times(\partial\Omega\setminus\Gamma_{\mathrm{c}})$. \autoref{lem:terminal_closure} gives the terminal inequalities on $\overline\Omega\setminus(\Gamma_{\mathrm{c}}\cup\overline{\Gamma_{\mathrm{abs}}})$. \autoref{prop:comparison} gives \textup{(C)}, and \autoref{thm:viscosity_limit}\textup{(iii)} gives the convergence and identifies $\overline u=\underline u$ on $Q_{\mathrm{c}}$ as the unique solution there.

\section{Numerical details}\label{app:numerical_provenance}

In this appendix, we collect the settings shared by all runs, then the parameters of each benchmark, a cost comparison with a five-point scheme, and the census of row types and fans.

The numerical experiments use the open-source MFGArchon framework \cite{MFGArchon2026} (post-v0.21) and the positive-Laplacian fan-upwind construction analyzed in the main text. Low-discrepancy clouds are generated from Sobol sequences \cite{Niederreiter1992}, inset from the boundary by a fraction of the nominal spacing. The base kernel $K$ is Gaussian and $h_{\mathrm{KDE}}$ is an absolute length in all studies. The source-free exact MFG uses the spectral realization, which is \eqref{eq:kde_operator} itself, at $h_{\mathrm{KDE}}=0.016$, chosen against $m^*$'s feature scale $1/\omega=0.0796$. The evacuation uses the ghost-image reconstruction $\mcalK^{\mathrm{gh}}_{h_{\mathrm{KDE}}}$ at $h_{\mathrm{KDE}}=0.49$, close to the nominal spacing of the $N_{\mathrm{int}}=1200$ cloud, held at that value on the finer clouds so that the refinement ladder varies the collocation cloud alone. In the floor-plan, source-free and evacuation runs, a row first takes the ball of radius $\delta$, about $2h$, and then the pool of radius $5\delta/2$, about $5h$, with $\kappa=3/5$. At this $\delta$, a pool of $2.25h$ refuses nodes far from $\Gamma_{\mathrm{c}}$. From $3h$, on the $31$ clouds we scanned, every refused node lies within about $2h$ of $\Gamma_{\mathrm{c}}$, and from $3h$ to $5h$ their number falls by about a quarter. A larger radius increases the error. Scaling $\delta$ and the pool together by $5/2$ multiplies the floor-plan error by $2$ to $3$ on $u_1$ and by $8$ to $20$ on $u_2$ from $N_{\mathrm{int}}=1300$ to $5200$. The transfer $\mcalI^{\mathrm{vec}}$ uses inverse-distance weights over the four nearest nodes on the reflecting square, and in the evacuation piecewise-linear interpolation on a Delaunay triangulation of the cloud, with the nearest node outside its hull. All randomized clouds, particle resampling, and Euler--Maruyama increments are seeded.

\subsection{The LQ benchmark}\label{app:rate_experiments}

The solver is the fan-upwind scheme with the exact fan policy \eqref{eq:howard_alpha}, the positive simplicial-fan advection $A_{\mathrm{fan}}^\alpha$ and the positive-Laplacian program, and it involves none of the KDE and particle terms. This domain has no reflecting wall, and a row takes nearest nodes in place of a pool. The Laplacian uses the $21$ nearest, enlarged to $31$, $41$ and $61$ where the program is infeasible, and the fan searches over the $40$ nearest. The $21$, $61$ and $40$ nearest nodes reach about $2.6h$, $4.4h$ and $3.6h$. We refine space and time jointly, with $\Delta t=h/40$. Errors are taken at $t=T/2$.

A jittered cloud places some nodes so close to the frame, within $0.3h$, that the two frame nodes bracketing such a node subtend more than $120^\circ$ with no node between them. Such rows make up $1.5$--$3.6\%$ of the interior on the jittered clouds of this subsection (none on the Sobol clouds). Their fans use the smallest aperture cap that admits a complete fan, at most $140^\circ$, which \autoref{def:fan_complete} allows, and then $C_\theta=2.92$. All other rows keep the aperture cap $120^\circ$. We report errors in $L^2$ and in $L^\infty$, the norm that the discrete comparison principle and \autoref{cor:residual_value_transfer} control.

\subsection{Comparison with a five-point scheme}\label{app:complexity_table}

\begin{figure}[htbp]
    \centering
    \includegraphics[width=0.42\textwidth,height=0.21\textheight]{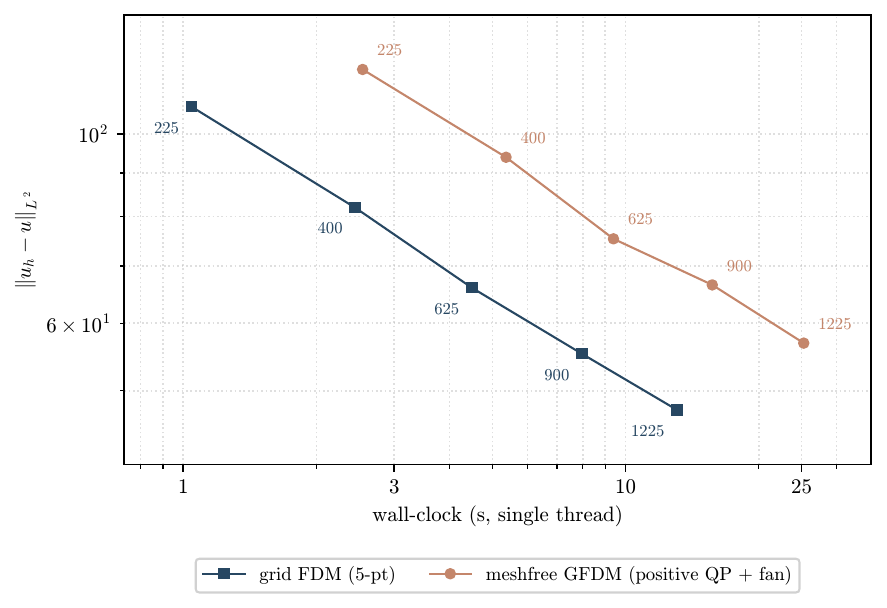}
    \caption{Error versus wall-clock for the positive-$\msfP_2$/guarded-fan HJB discretization and a five-point FDM on the decoupled 2D LQ problem. Five matched-DOF levels are shown from a single-thread run. Labels give interior DOFs. Particle and KDE costs are excluded.}
    \label{fig:fdm_pareto}
\end{figure}

On a single-thread run of the decoupled LQ problem (Sobol seed $42$, five levels, $\Delta t=h/67$), both methods perform one sparse LU factorization per policy update, stop Howard's iteration on the same test, and count stencil assembly and solve. The offline positive-QP and fan construction is inside the meshfree time, and the FDM has no counterpart to it. At $N_{\mathrm{int}}=1225$ the meshfree method spends $7.52$~s offline against $0.01$~s, and $17.76$~s in the solve against $13.06$~s. Across the five levels the offline share falls from $55.9\%$ to $29.7\%$ while the solve-only ratio rises from $1.08$ to $1.36$, so the total ratio falls with refinement over this range. The positive-$\msfP_2$ Laplacian stores $21$ nonzeros per interior row against the five-point stencil's five. The fitted $L^2$ EOCs are $0.867$ for the meshfree scheme and $1.010$ for the FDM (\autoref{fig:fdm_pareto}).

\subsection{Source-free exact MFG}\label{app:source_free_benchmark}

Each run activates the HJB solve, drift interpolation, reflected Euler--Maruyama particles, density reconstruction, and twelve fictitious-play steps with $\eta_k=1/(k+1)$. It uses $N_p=10^6$ and, for this reflecting square, the spectral realization on a $512^2$ transform grid. On samples of $m^*$, the $512^2$ and $1024^2$ transforms agree to four digits in the reported norm, so the reconstruction is grid-converged. The nominal spacing is $h=N_{\mathrm{int}}^{-1/2}$ and $\Delta t=T/N_t$. Within each seed, the initial ensemble and the increment stream are seeded identically across levels. For seeds $42$, $1$, and $2$, define the mid-time interior-quadrature error
\begin{equation*}
    E_m^{(s)}:=\norm{m_h^{(s)}(T/2)-m^*}_{L^1(\Omega)}.
\end{equation*}
The levels $N_{\mathrm{int}}=400,900,1600,2500$ use $N_t=125$, and the sample standard deviation over seeds $42$, $1$ and $2$ is at most $3.9\%$ of the mean.

Since $\norm{m^*}_1=1$, the finest value is a $3.56\%$ relative $L^1$ error. The time step is held because refining it with the cloud does not isolate the cloud. At $N_t=20,45,80,125$ on these levels the mean error falls from $0.136$ to $0.0356$, and particles driven by the exact drift at the same time steps give $0.143,\ 0.073,\ 0.045,\ 0.032$ on seed $42$.

\subsection{Evacuation discretization}\label{app:evac_setup}

\autoref{tab:evac_discretization} collects the parameters at $N_{\mathrm{int}}=1200$. The nine runs use five seeds at $N_{\mathrm{int}}=1200$, two of them also at $1800$ and $2400$, and all particles have exited by $T=40$. An exact segment test removes every edge that meets solid material at assembly.

\emph{Wall reflection.} One Euler step can span the wall, since $\sigma\sqrt{\Delta t}+\abs{\nabla u}\Delta t$ is comparable to $w$. Contact is therefore solved in closed form, on the two axis-aligned boxes of the door wall and the quadratic of the obstacle, and the reflection of \eqref{eq:euler_maruyama} iterates up to a bounce cap. On all nine clouds no segment of any reflection polyline meets wall material and the cap is never reached.

\begin{table}[!htbp]
    \footnotesize
    \centering
    \caption{Discretization parameters for the evacuation at $N_{\mathrm{int}}=1200$.}
    \label{tab:evac_discretization}
    \begin{tabular}{llc}
        \toprule
        Stage & Parameter & Value \\
        \midrule
        Interior cloud & Sobol quasi-Monte Carlo points $N_{\mathrm{int}}$ & $1200$ \\
        & further nodes in $[12.25,15.75]\times[4.25,7.75]$ & $280$ \\
        & outer-boundary sampling margin & $0.15$ \\
        \addlinespace
        Boundary cloud & exit nodes at spacing $0.1$ & $80$ \\
        \addlinespace
        Stencil & pool radius $\delta_i$, with $\delta=0.911$ at $N_{\mathrm{int}}=1200$ & $5\delta/2$ \\
        & reach, so that $\kappa=3/5$ & $3\delta/2$ \\
        \addlinespace
        Outer MFG & outer iterations & $6$ \\
        & fictitious-play step $\eta_k$ & $1/(k{+}1)$ \\
        \addlinespace
        Forward FP & particles $N_p$ & $2 \times 10^4$ \\
        & time steps $N_t$ & $162$ \\
        & horizon $T$ & $40$ \\
        & diffusion $\sigma$ & $0.5$ \\
        \bottomrule
    \end{tabular}
\end{table}

\subsection{Row and fan census}\label{app:row_census}

\autoref{tab:row_census} gives the outcomes of the row classification on the evacuation clouds. \autoref{tab:fan_census} summarizes the rest of the census: the installed fans on the LQ benchmark clouds, and the apertures and obliqueness margins on the evacuation clouds. Each installed fan is found by a search over the visible pool, which keeps adjacent rays at least $25^\circ$ apart and caps the aperture. Off-center Laplacian weights go down to $-1.3\times10^{-8}$, at solver tolerance. Setting them to zero moves each off-diagonal entry of $A^\alpha$ by at most $\tfrac{\sigma^2}{2}\cdot1.3\times10^{-8}$ and gives rows with the signs that \autoref{thm:constructed_operator_monotone} assumes. The runs use the weights as computed.

\begin{table}[htbp]
    \centering
    \footnotesize
    \caption{Outcomes of \autoref{def:row_classification} over the nine-cloud evacuation ensemble. Rows of type \textup{(P)} are listed by stencil, initial or visibility-safe enlarged.}
    \label{tab:row_census}
    \begin{tabular}{llr}
        \toprule
        Outcome & Stencil & Count \\
        \midrule
        \textup{(P)} PDE row & initial & $14{,}939$ \\
        \textup{(P)} PDE row & visibility-safe enlarged & $256$ \\
        \textup{(R)} mirror row & & $1684$ \\
        \textup{(D)} absorbing-Dirichlet identity row & & $0$ \\
        Copy row next to a boundary vertex & & $41$ \\
        \midrule
        Interior-cloud rows & & $16{,}920$ \\
        \bottomrule
    \end{tabular}
\end{table}

\begin{table}[htbp]
    \centering
    \footnotesize
    \caption{Fan and obliqueness census. Upper block: the five Sobol levels of the LQ benchmark ($3375$ interior rows). Lower block: the nine-cloud evacuation ensemble ($16{,}920$ interior rows), with signs stated to solver tolerance.}
    \label{tab:fan_census}
    \begin{tabular}{lp{0.48\textwidth}}
        \toprule
        Quantity & Value \\
        \midrule
        \multicolumn{2}{l}{\emph{LQ benchmark clouds}} \\
        Adjacent-ray angle & $25.0^\circ$ to $119.7^\circ$ \\
        Fan constant $C_\theta$ & $1.99$ \\
        \midrule
        \multicolumn{2}{l}{\emph{Evacuation ensemble}} \\
        Largest installed aperture & at most $120^\circ$ on PDE rows, $134.8^\circ$ on mirror rows \\
        $\min\widehat\gamma_{i,h}$ at $\gamma_0=1.25$ & $1.25$, set by three clouds, and from $1.28$ to $1.56$ on the six where the constraint is slack \\
        Moment residual & below $1.4\times10^{-11}$ \\
        Smallest off-center weight & $-1.3\times10^{-8}$ \\
        \bottomrule
    \end{tabular}
\end{table}

\end{document}